\documentclass{article}
\usepackage{amsfonts}
\usepackage{amssymb}
\usepackage{bbm}
\usepackage{comment}
\usepackage{authblk}
\usepackage{enumitem}
\usepackage{academicons}
\usepackage{dsfont}
\usepackage{amsthm}
\usepackage{amsmath}
\usepackage{commath}
\usepackage{url}
\usepackage{float}
\usepackage{epsfig}
\usepackage{subfigure}
\usepackage{adjustbox}
\usepackage{mathrsfs}  
\usepackage[ruled,vlined,linesnumbered]{algorithm2e}
\usepackage{bm}
\usepackage{graphicx}
\usepackage{a4wide}
\usepackage{todonotes}
\usepackage{mathtools}
\makeatletter
\@namedef{subjclassname@2020}{%
	\textup{2020} Mathematics Subject Classification}
    \newcommand{\subjclass}[2][2020]{%
  \let\@oldtitle\@title%
  \gdef\@title{\@oldtitle\footnotetext{#1 \emph{Mathematics subject classification (MSC 2020):} #2}}%
}
\makeatother

\usepackage[T1]{fontenc}

\DeclareMathOperator{\N}{\mathbb{N}}

\DeclareMathOperator{\Z}{\mathbb{Z}}
\newtheorem{thm}{Theorem}
\newtheorem*{conj}{Conjecture}
\newtheorem{lem}[thm]{Lemma}
\newtheorem{defi}{Definition}

\newtheorem{prop}[thm]{Proposition}
\newtheorem{cor}[thm]{Corollary}
\newtheorem{lemma}[thm]{Lemma}

\theoremstyle{definition}

\DeclareMathOperator{\vol}{Vol}

\newcommand{\ls}{\leqslant}
\newcommand{\gs}{\geqslant}

\newcommand{\orcid}[1]{\href{https://orcid.org/#1}{\textcolor[HTML]{A6CE39}{\aiOrcid}}}

\usepackage{graphicx} % Required for inserting images

\title{The bad and rough rotation is Poissonian}
\author{Manuel Hauke-Treuer}

\date{\today}

\begin{document}

\newcommand{\Addresses}{{% additional braces for segregating \footnotesize
  \bigskip
  \footnotesize

  M.~Hauke, \textsc{Department of Mathematical Sciences, NTNU Trondheim},\par 
\textsc{Sentralbygg 2, 7034 Trondheim, Norway}\par\nopagebreak
  \textit{E-mail address} \texttt{hauke@math.tugraz.at}\par\nopagebreak
  \textit{Orcid:} \texttt{0000-0002-6244-0285}
  }}

\maketitle

\abstract{
%Motivated by the Berry--Tabor Conjecture and the seminal work
%of Rudnick--Sarnak, 
The fine-scale properties of sequences $(a_n\alpha)_{n \in \mathbb{N}} \mod 1$ with $(a_n)_{n \in \mathbb{N}} \subseteq \mathbb{N} $ and $\alpha$ irrational have been extensively studied in the last decades. In this article, we prove that for $(a_n)_{n \in \mathbb{N}}$ arising from the set of rough numbers with explicit roughness parameters and any badly approximable $\alpha$, $(a_n\alpha)_{n \in \mathbb{N}} \mod 1$ has Poissonian correlations of all orders, and consequently, Poissonian gaps. This is the first known explicit sequence $(a_n\alpha)_{n \in \mathbb{N}} \mod 1$ with any of these properties. Further, we show that the statement is false for Lebesgue almost every $\alpha$, thereby disproving a conjecture of Larcher and Stockinger [Math. Proc. Camb. Phil. Soc. 2020].
The method of proof makes use of an equidistribution result in $\mathbb{Z}_d$ in Diophantine Bohr sets which is established by appealing to a novel method of examining random walks on Ostrowski digit expansions.
}

\section{Introduction and statement of results}

A sequence $(x_n)_{n \in \mathbb{N}} \subseteq [0,1)$ is said to be uniformly distributed modulo 1 if asymptotically, each interval $I$ contains its fair share of points, i.e. for all $I = [a,b] \subseteq [0,1)$, we have

\[\lim_{N \to \infty}\frac{1}{N}\#\{n \leq N: x_n \in I\} = \lambda(I),\]
where $\lambda$ denotes the Lebesgue measure. This notion goes back at least to Weyl’s seminal paper \cite{Weyl} and was a very intensively studied concept in the last century. We remark that the property for a sequence to be uniformly distributed can be seen as a measure of (pseudo-)randomness: Choosing a sequence of i.i.d. random variables $(X_n)_{n \in \mathbb{N}}$ with $X_1 \stackrel{d}{=} U([0,1))$, we obtain that 
$(X_n)_{n \in \mathbb{N}}$ is almost surely uniformly distributed.

Another, stronger (see \cite{ALP,GL,HZ_higher_order,Mark_manifold,Steinerberger}) pseudo-random property  for a sequence is having \textit{Poissonian correlations of $k$-th order}: A sequence is said to have Poissonian correlations of $k$-th order for $k \geq 2$ if for any axis-parallel rectangle $\bm{R} \subseteq \mathbb{R}^{k-1}$, we have\footnote{There are other, equivalent ways to define Poissonian correlations e.g. by taking smooth test function instead of rectangles or considering the tuples $(x_{n_1}- x_{n_2},x_{n_2}- x_{n_3},\ldots x_{n_{k-1}}- x_{n_k})$ instead. For a proof of the equivalences between these definitions, we refer the reader to \cite[Appendix A]{HZ_higher_order}.}
\[R_k(\bm{R},N) := 
\frac{1}{N}\#\left\{\bm{n} \in [1,N]^k_{\neq}: (x_{n_1}- x_{n_2},x_{n_1}- x_{n_3},\ldots, x_{n_1}- x_{n_k}) \in \frac{\bm{R}}{N}\right\} \underset{N \to \infty}{\to} \vol(\bm{R}),
\]
where $[1,N]^k_{\neq}$ stands for all $k$-tuples with pairwise distinct entries in $[1,N]$, and $\vol$ denotes the $(k-1)$-dimensional Lebesgue measure. The correlation functions $R_k$ can be defined in the analogous way for triangular arrays $((x_n)_{n \leq N})_{N \in \mathbb{N}}$ as well, see e.g. \cite{Mark_manifold} for such a generalization.
While a sequence of i.i.d. random variables as above has almost surely Poissonian correlations of $k$-th order (see e.g. \cite[Appendix B]{AHZ}), there are only very few \textit{explicit} sequences known to have these properties: 
 El-Baz, Marklof and Vinogradov \cite{El-baz} showed that $\{\sqrt{n}\}_{n \in \mathbb{N}\setminus\{m^2: m \in \mathbb{N}\}}$ (here and throughout this article, $\{\cdot\}$ denotes the fractional part of a real number) has Poissonian pair correlations (PPC), i.e., Poissonian correlations of second order.
Lutsko, Sourmelidis and Technau \cite{LST} showed that $\{\alpha n^{\theta}\}$ has PPC for all $\alpha$ and all $\theta < 14/41$, a result improved by Radziwi\l\l \;and Shubin \cite{Rad_shu} to all $\theta < 43/117$. For higher-order correlations, even less is known: Lutsko and Technau \cite{LT_n_theta}
showed that $\{\alpha n^{\theta}\}_{n \in \mathbb{N}}$ has Poissonian correlations of $k$-th order for all $2 \leq k \leq m$ when $\theta < \theta_m$ where $\theta_m \to 0$ as $m \to \infty$, and $\{\alpha(\log n)^{A}\}_{n \in \mathbb{N}}, A > 1$ has Poissonian correlations of $k$-th order for all $k \geq 2$ \cite{LT_log}.\\

The main motivation for studying (Poissonian) correlations is to determine the gap distribution of a given sequence:
Writing $x_1^{(N)} \leq x_2^{(N)} \ldots \leq x_N^{(N)}$ for the ordered sample of $\{x_n: 1 \leq n \leq N\}$, we say that a sequence $(x_n)_{n \in \mathbb{N}} \subseteq [0,1)$ has  Poissonian gap distribution (or sometimes also called Exponential gaps) if for all $[a,b] \subset [0,\infty)$, we have
\[\lim_{N \to \infty} \frac{1}{N}\#\left\{1 \leq n \leq N-1: x_{n+1}^{(N)} - x_{n}^{(N)} \in \frac{[a,b]}{N}\right\} = \int_{a}^b e^{-t} \,\mathrm{d}t.\]
Analyzing the gap distribution directly is notoriously difficult.
While it is known that having Poissonian gaps does not imply uniform distribution and thus also not Poissonian correlations of any order \cite{AHZ}, the most promising approach to showing Poissonian gaps usually goes by proving that a sequence has Poissonian correlations of $k$-th order for all $k\geq 2$, which then implies that a sequence has Poissonian gaps, see \cite[Appendix A]{Kur_Rud} for a proof of this implication.\\

Finding the gap distribution of a sequence is considered generally challenging:
%for some mathematicians in this area to be ``the holy grail'' \cite[p.3]{LST}: 
There are various deep and widely open conjectures in this area that include the Berry--Tabor Conjecture from theoretical physics (see \cite{Berry_Tabor, Marklof_survey,Rudnick_survey} for more details), and the Rudnick--Sarnak--(Zaharescu) Conjecture \cite{rud_sar,rhz} that was repeated by Heath--Brown \cite{hb}: 
Under a Diophantine condition on $\alpha$ (which is $\alpha$ being badly approximable in \cite{rud_sar}), $\{n^2\alpha\}_{n \in \mathbb{N}}$ has Poissonian gaps, and clearly, a bound on the Diophantine exponent of $\alpha$ is necessary.
%\footnote{Note that for the case $d = 1$ for trivial reasons, $\{n^d\alpha\}_{n \in \mathbb{N}}$ does neither have Poissonian correlations of any order nor Poissonian gaps, regardless of the value of $\alpha$.}. 
While there has been partial progress \cite{hb,lut,rhz,tw,truelsen}, the conjecture itself remains widely open.

As to be expected, only a few examples are known where the gap distribution of a sequence in $[0,1)$  has been determined: Lutsko and Technau showed Poissonian gaps for $\{\alpha(\log n)^{A}\}_{n \in \mathbb{N}}, A > 1, \alpha \in \mathbb{R}$ \cite{LT_log}. Other sequences where a nontrivial, but non-Poissonian gap distribution was determined include notably $\{\sqrt{n}\}_{n \in \mathbb{N}}$ by Elkies and McMullen \cite{elk_mc}, and $(\{\log_b n\})_{n \in \mathbb{N}}$ by Marklof and Strömbergsson \cite{marklof_log}.\\

Motivated especially by the Rudnick--Sarnak--Zaharescu Conjecture, various sequences of the form $\{a_n\alpha\}_{n \in \mathbb{N}}$, where $(a_n)_{n \in \mathbb{N}}$ is an increasing integer sequence and $\alpha$ is an irrational number, have been studied. Since finding explicit examples is considered very challenging, this happened so far only from a metric point of view:
 Poissonian correlations respectively gaps have been proven for $\{a_n\alpha\}_{n \in \mathbb{N}}$ for Lebesgue-almost every $\alpha$ for various integer sequences $(a_n)_{n \in \mathbb{N}}$. While Rudnick and Zaharescu  \cite{RZ} showed Poissonian gap distributions for almost all $\alpha$ whenever $(a_n)_{n \in \mathbb{N}}$ is lacunary, most other results established so far consider only the most tractable case $k = 2$, i.e. Poissonian pair correlations: Rudnick and Sarnak proved that $\{n^d\alpha\}_{n \in \mathbb{N}}, d \geq 2$ has PPC for almost all $\alpha$ \cite{rud_sar}. A more general examination of the metric theory of $\{a_n\alpha\}_{n \in \mathbb{N}}$ was given by Aistleitner, Lewko and Larcher \cite{all} who provided a sufficient criterion on the additive energy of $(a_n)_{n \in \mathbb{N}}$ in order for $\{a_n\alpha\}_{n \in \mathbb{N}}$ to have almost surely PPC. In the appendix of the same article, Bourgain provided a necessary condition. For more results in this direction, we refer to e.g. \cite{alt,chow_gafni_walker,bloom_walker,Lach_Tec19,felipe}. Finally, we mention the work of Walker \cite{walker_primes} who showed that the prime rotation $\{p_n\alpha\}_{n \in \mathbb{N}}$ satisfies PPC for Lebesgue-almost no $\alpha$.
We remark that all these statements are of a metrical nature, or make assumptions on $\alpha$ where no explicit $\alpha$ is known to satisfy the assumed properties.
To the best of the author's knowledge, prior to this article no explicit sequence of the form $\{a_n\alpha\}_{n\in \mathbb{N}}$ with $(a_n)_{n \in \mathbb{N}} \subseteq \mathbb{N}$ was known to have Poissonian correlations even for $k =2$, let alone Poissonian gaps.\\

In this article, the set of irrational numbers $\alpha$ considered is the set of badly approximable numbers, and thus coincides with the assumption of the Rudnick--Sarnak Conjecture in \cite{rud_sar}. An irrational $\alpha$ is called badly approximable whenever $\liminf_{n \to \infty}n \lVert n \alpha\rVert > 0$ (where $\lVert.\rVert$ stands for the distance to the nearest integer) or equivalently, the partial quotients of the continued fraction expansion of $\alpha$ are uniformly bounded. The set of badly approximable numbers has Lebesgue measure $0$, but Hausdorff dimension $1$, and contains prominent explicit examples such as the Golden Ratio $\Phi= \frac{1 + \sqrt{5}}{2}, \sqrt{2}$, or in general, all irrational solutions arising from a quadratic polynomial with integer coefficients.
The assumption of $\alpha$ being badly approximable will be essential in our arguments and will be shown to be almost necessary.\\

The sequence $(a_n)_{n \in \mathbb{N}}$ considered in this article arises from the classical generalization of prime numbers to rough numbers. We say that an integer $n$ is $z$-rough if the smallest prime divisor of $n$, denoted by $P^-(n)$, is larger than $z$. While for $\sqrt{n} < z < n$, this is equivalent to saying that $n$ is prime, the generalization to smaller $z$ often allows to study denser sets of numbers that are ``almost prime'' in the sense that they have very few prime divisors, but the occurrence of such numbers is easier to detect or quantify, especially with sieve-theoretic considerations. Many famously open number-theoretic questions, including the Twin Prime conjecture \cite{gpy}, the Goldbach Conjecture \cite{chen_gold,renyi},  and Diophantine approximation with prime denominator \cite{harm_primes,irv}, have been (partially) solved in the easier case when considering suitably rough numbers instead of primes themselves. This general philosophy will be continued in this article, albeit we remark that having some \textit{lower} bound on the roughness condition is also necessary (mainly to ensure a certain growth rate on $(a_n)_{n \in \mathbb{N}}$) to have a chance for $\{a_n\alpha\}_{n \in \mathbb{N}}$ to have any Poissonian property -- after all, it is immediate to see that for $a_n = n$, $\{a_n\alpha\}_{n \in \mathbb{N}}$ does not have any Poissonian property, regardless of the choice of $\alpha$.
We define the sequence of rough numbers with roughness function $f$ in the following way:

\begin{defi}
For given $f: \N \to [0,\infty)$, we define
\[\mathcal{R}_f := \{n \in \N: P^{-}(n) > f(n)\}\]
and let $(a_n^{(f)})_{n \in \N}$ be $\mathcal{R}_f$ written in increasing terms.\\
Further, we define the triangular array
$(a_{n,x}^{(f)})_{n \leq x}$ of (uniformly) rough numbers by 
\[\{a_{n,x}^{(f)}: 1 \leq n \leq x\} := \{1 \leq n \leq x: P^{-}(n) > f(x)\}.\]
We define the function that counts (uniformly) rough numbers by 
\[\Phi(x,z) := \#\{1 \leq n \leq x: P^{-}(n) > z\}.\]
\end{defi}
\vspace{5mm}

With this definition at hand, we state the main result of this article.
\vspace{3mm}

\begin{thm}\label{main_thm}
Let $f(x) = x^{r(x)}$ where 
$r(x)$ decreases monotonically to $0$, but satisfies 
$r(x) \gg_A \frac{(\log \log x)^A}{\log x}$ for all $A > 0$.
Then for any badly approximable $\alpha$, the sequence 
$(a_n^{(f)}\alpha)_{n \in \N} \pmod 1$ has Poissonian correlations of all orders. Consequently, the sequence 
$(a_n^{(f)}\alpha)_{n \in \N} \pmod 1$ has Poissonian gaps.\\
Furthermore, under the same assumptions on $f$, the triangular array of dilates of rough numbers 
$(a_{n,x}^{(f)}\alpha)_{n \leq x} \pmod 1$ has Poissonian correlations of all orders, and thus also Poissonian gaps.\\
\end{thm}

As concrete example, we might choose in Theorem \ref{main_thm} those $n$ that are $\exp(\sqrt{\log n})$-rough, and take $\alpha = \frac{1 + \sqrt{5}}{2}$, providing an explicit example of such a sequence.
Further, we note that knowing that all correlations are Poissonian, one can also recover that the $k$-th neighbour spacings (i.e. $x_{n+k}^{(N)} - x_{n}^{(N)},\, 1 \leq n \leq N$) are also Poissonian, for every fixed $k\geq 1$. For a proof of this fact, we refer once more to \cite[Appendix A]{Kur_Rud}.\\

As a converse result to Theorem \ref{main_thm}, we show that the Diophantine property of $\alpha$ being badly approximable is not an artifact of the proof: If $\alpha$ is ``slightly better approximable'', then the statement becomes false even for the simplest case $k = 2$:

\begin{thm}\label{dioph_ass_thm}
    Let $f(x) = x^{r(x)}$ where $r(x)$ decreases monotonically to $0$ and let $\alpha$ be an irrational number satisfying
    \begin{equation}\label{dioph_ass}\liminf_{n \to \infty} n \lVert n\alpha\rVert \log (f(n)) = 0.\end{equation}
    Then $(a_n^{(f)}\alpha)_{n \in \mathbb{N}} \pmod 1$ does not have Poissonian pair correlations.\\ In particular, for Lebesgue-almost every $\alpha$, the result of Theorem \ref{main_thm} fails.
\end{thm}

The proof of Theorem \ref{dioph_ass_thm} is a straightforward adaptation of Walker’s argument for primes \cite{walker_primes}. We remark that it implies that the fine-scale statistics of $\{a_n^{(f)}\alpha\}_{n \in \mathbb{N}}$ indeed rely delicately on the Diophantine properties of $\alpha$. This is in stark contrast to the real-valued sequences $((\log n)^A)_{n \in \mathbb{N}}, (n^{\theta})_{n \in \mathbb{N}}$ examined in \cite{LST,LT_n_theta,LT_log,Rad_shu} where the Diophantine properties of $\alpha$ (as maybe expected for real-valued sequences $(a_n)_{n \in \mathbb{N}}$) do not play a role whatsoever, and can be chosen arbitrarily (even rational).\\

The combination of Theorems \ref{main_thm} and \ref{dioph_ass_thm} further provides a warning remark to be very cautious in deducing that a result known for almost all $\alpha$ (such as the flurry of results mentioned above) implies anything for the Lebesgue null set of badly approximable numbers. Indeed, Larcher and Stockinger \cite{LS20} conjectured the following:

\begin{conj}[Larcher, Stockinger, 2018]
Let $(a_n)_{n \in \mathbb{N}}$ be an integer sequence.
If for (Lebesgue) almost all $\alpha$ the pair correlations of $\{a_n \alpha\}_{n \in \mathbb{N}}$ are not Poissonian, then the pair correlations of this sequence are not Poissonian for any $\alpha$.
\end{conj}

A combination of Theorem \ref{main_thm} and (the ``in particular'' part of) Theorem \ref{dioph_ass_thm} implies immediately that the Conjecture turns out to be false, witnessed by $(a_n^{(f)})_{n \in \mathbb{N}}$ with $f$ as in Theorem \ref{main_thm}.\\

\subsubsection*{Remarks on the theorems, primes in short intervals, and related works}

We point out that the main novelty of Theorem \ref{main_thm} compared to earlier works on $(a_n\alpha)_{n \in \mathbb{N}}$ is the \textit{explicit} construction of the sequence $(a_n \alpha)_{n \in \mathbb{N}}$. Furthermore, the choice of various functions for the roughness allows us to take a family of different sequences $(a_n)_n$, while the growth rate is provided by $a_n/n \sim r(n) \log n$, which allows for various growth rates, as long as
$(\log \log n)^A \ll_A \frac{a_n}{n} = o(\log n)$.
 %Concerning the assumptions on $f$ made in Theorem \ref{main_thm}, we make the following remarks: 
 We believe the lower bound on $r$ in Theorem \ref{main_thm} to be sharp: This follows from the proof (see Section \ref{sec_heur} for a proof sketch) where averages of multiplicative functions need to be considered (this is not due to the proof method, but equivalent to having Poissonian correlations). For the case $k = 2$, the necessary condition \begin{equation}\label{avg_rmk}\frac{1}{\#\mathcal{B}}\sum_{n \in \mathcal{B}}\frac{n}{\varphi(n)} \ll 1\end{equation} 
 (with $\varphi$ the Euler totient function) arises where $\mathcal{B} = S(\alpha,x) \subseteq [1,x]$ is a so-called Bohr set of cardinality $\approx r(x) \log x$. Assuming the (reasonable) existence of badly approximable $\alpha$ with infinitely many $n \in \mathbb{N}$ such that $n \in \bigcup_{x \in \mathbb{N}}S(\alpha,x)$ with $ \frac{n}{\varphi(n)} \gg \log \log n$, such an averaging result becomes impossible when $r(x) \log x = o(\log \log x)$, since even one such exceptional $n \in \mathcal{B}$ would spoil \eqref{avg_rmk}. Similarly, correlations of higher order prevent us from considering $r(x) \log x < (\log \log x)^A$ for any fixed $A > 0$.\\
 
On the other hand, the assumption of $r(x) \to 0$ arises from a classical sieve-theoretic obstruction: if we knew a uniform version of the Hardy--Littlewood $k$-tuple Conjecture \cite{HL}, that contains e.g. the Twin Prime Conjecture as a special case, the result would also be true without $r(x) \to 0$. This leads us to the following Conjecture:

\begin{conj}[Poissonian Prime rotation]\label{conj_primes}
    Let $\alpha$ be a badly approximable number. Then the prime rotation sequence $(p_n\alpha)_{n \in \mathbb{N}} \mod 1$ has Poissonian correlations of all orders, and consequently, has Poissonian gap distribution.
\end{conj}

We remark that the convergence towards the Poissonian gap distribution appears numerically to be very slow, both for the rough numbers covered by Theorem~\ref{main_thm} and in the conjectured prime case. At computationally feasible finite cutoffs, the expected asymptotic behaviour is still strongly obscured by effects arising from the Three Gap Theorem, making the gap distribution look distinctly non-Poissonian at first sight. This is largely due to the comparatively high density of the primes and rough numbers:
for instance, for the much sparser sequence of squares, a Poissonian-looking gap distribution is already visible at substantially smaller cutoffs. For a more detailed explanation of this finite-scale phenomenon, together with further numerical experiments supporting the Poissonian Prime rotation Conjecture, we refer the reader to the Appendix.\\

Note that there is a genuine obstruction for proving the Conjecture by our method with currently known estimates: We employ asymptotic sieve estimates to differences of rough numbers. In the case of prime numbers, these estimates could only be established if we knew a uniform version of the Hardy--Littlewood $k$-tuple Conjecture \cite{HL}, that contains e.g. the asymptotic Twin Prime Conjecture as a special case. 
% We remark at this point that Freiberg, Kurlberg and Rosenzweig studied Poissonian gaps of the sum of two squares by conjectural on an analogue of the Hardy--Littlewood $k$-tuple Conjecture in the corresponding setup \cite{Kur_conj}.\\
In order to circumvent such widely open conjectures, one would (already for $k = 2$) need to count the number of prime pairs $p,p' \leq x$ with $p-p' = n$ on a weighted average over the Bohr set
$\mathcal{B}(\alpha,x) = \{n \leq x: \lVert n\alpha \rVert \leq \tfrac{\log x}{x}\}$, see Section \ref{sec_heur} for context.
In the method employed for Theorem \ref{main_thm}, we apply these estimates point-wise (in $n$), not making use of the averaging information at this stage of the proof. We can only establish this point-wise bound because of the roughness assumption $f(x) = x^{o(1)}$ instead of primes (i.e. $f(x) = x^{1/2}$).
Nevertheless, the Hardy--Littlewood $k$-tuple Conjecture is true for almost all integers (see \cite{aveng,mikawa}), and thus there might be some hope to adapt the averaging techniques to also work on the corresponding Bohr sets.
While Diophantine Bohr sets have some rigid additive structure (which is exploited very heavily in the proof of Theorem \ref{main_thm} at a later step when averaging over the arising singular series) that might help to establish results ``on average'', the exponential sparsity $\#\mathcal{B}(\alpha,x) \asymp {\log x}$ of the Bohr sets is way too thin for even the most modern tools such as \cite{aveng}.\\

\textbf{Poissonian correlations on the real line -- related works.} We note that so far, the focus of $\{a_n\alpha\}_n$ concentrated, even from a metric perspective, on either monomial sequences, or arbitrary sequences, with the exception being Walker's article \cite{walker_primes} on primes. However, various sequences of number-theoretic interest have been studied on the real line:
The notion of Poissonian correlations can also be transferred to the real line, by normalizing an increasing sequence $(x_n)_{n \in \mathbb{N}} \subseteq \mathbb{R}$ at each step $N$ by the $x_N-x_1$; e.g., one can study the primes, asking about the distribution of $(p_{n+1}-p_n)/\log N$ for $1 \leq n \leq N$. Gallagher \cite{Gall_primes} showed that under a uniform version of the Hardy--Littlewood Conjecture (which coincides with the version we would need to assume to prove the Poissonian Prime rotation Conjecture) that the primes are Poissonian (as objects on the real line, as described above). 
Notably, Montgomery and Soundararajan \cite{mont_sound} conjectured that on slightly longer ranges (i.e., counting correlations at scale $H \geq (\log N)^{1+\delta}$), the primes are expected to converge to a Gaussian distribution, with variance $H \log (N/H)$ -- see also recent advancements e.g. \cite{vivian,sunkai} on related matters.
Similarly, Freiberg, Kurlberg and Rosenzweig \cite{Kur_conj} showed Poissonian behaviour (on the real line) for sums of two squares, subject to an analogous Hardy--Littlewood-type Conjecture on the sums of two squares, mirroring the result of Gallagher \cite{Gall_primes}.\\

We remark that the behaviour of Poissonian correlations on the real line is somewhat different to studying \textit{sequences} in $[0,1]$: First of all, the simple observation that the sequence is increasing implies that previous existing gaps are not split, but only marginally rescaled, while in the case of a fixed sequence in $[0,1)$, an earlier gap might be split at any time. Clearly, the re-normalizing factor allows these sequences to be interpreted as triangular arrays in $[0,1)$, a more general setup in which further results of number-theoretic interest are known, such as the works of Hooley \cite{HooleyI,hooleyII,hooleyIII} on $\mathbb{Z}_q^{*}$ who showed Poissonian behaviour for $q/\varphi(q) \to \infty$, as well as the works of Kurlberg and Rudnick \cite{Kur_Rud,Kur_nsqf} on quadratic residues for highly-composite $q$.

% We note however that the proof method presented for Theorem 1 allows us to obtain Conjecture (give it a number) conditionally under the same uniform version of the HL-Conj than Gallagher used in (...).

% Kurlberg and Rudnick showed Poissonian gaps for the quadratic residues of a highly composite modulus \cite{Kur_nsqf,Kur_Rud}; 

\subsection{Equidistribution mod $d$}

In order to estimate the Poissonian correlations for rough numbers, after applying a sieve estimate (see Section \ref{sec_heur} for a proof layout), it becomes evident that one has to control the average of the well-known singular series

\[{\mathfrak S}(\bm{h}) =  {\mathfrak S}_{k-1}(\bm{h})= \prod_{p \leq f(x)}\left(1 - \frac{g_{\bm{h}}(p)}{p}\right)\]
where 
$g_{\bm{h}}(p) = \#\{0 \pmod p,h_{1} \pmod p,\ldots,h_{k-1} \pmod p\}$, with the average taken over the Bohr set
$(\mathcal{B}(\alpha,x))^{k-1}$ where
\[\mathcal{B}(\alpha,x) := \{1 \leq n \leq x: \{n\alpha\} \in [0,\rho_x]\},\]
with $x \rho_x \to \infty$ at a very slow rate. This can be established by understanding the divisor properties within $\mathcal{B}(\alpha,x)$ very well in order to obtain an asymptotics for the averaging described above. This is provided by the following statement.

\begin{thm}[Bohr set divisibility]\label{thm_bohr_div}
Let $\alpha$ be badly approximable and let $(I_x)_{x \in \N}$ be a sequence of intervals in $[0,1)$ of the form $[0,\rho_x]$ or $[1-\rho_x,1]$ with
$\lim_{x \to \infty} x\rho_x = \infty$,
and let $\mathcal{B}(x,I_x)$ denote the Bohr set
\[\mathcal{B}(x,I_x) := \mathcal{B}_{\alpha}(x,I_x) := \{1 \leq n \leq x: \{n\alpha\} \in I_x\}.\]
Then for every $d\in \N$ and every $a \in \mathbb{Z}_d := \mathbb{Z}/d\mathbb{Z}$, we have
\begin{equation}\label{bohr_div_asymp}
   \lim_{x \to \infty} \frac{\#\{n \in \mathcal{B}(x,I_x): n \equiv a \pmod d\}}{\#\mathcal{B}(x,I_x)} = \frac{1}{d}.
\end{equation}
\end{thm}

We remark that Theorem \ref{thm_bohr_div} (respectively the way it is applied later) might be of independent interest: Since Theorem \ref{thm_bohr_div} is very flexible (note that no speed of $x \rho_x \to \infty$ is necessary), it allows for averaging various multiplicative functions over Bohr sets. Thus, general de-correlating asymptotics such as

\begin{equation}\label{decor}\sum_{n \leq x} g(n)F_x(\{n\alpha\})
\sim \left(\frac{1}{x}\sum_{m \leq x} g(m)\right)\sum_{n \leq x} F_x(\{n\alpha\}),\quad x \to \infty,
\end{equation}
where the mass of $F_x$ is concentrated around $0$, and $g$ is a reasonably well-behaving multiplicative function, might become attainable for arbitrary badly approximable $\alpha$. We deliberately decided against stating a more general technical statement that provides sufficient criteria for the functions $g$ and $F_x$, such that \eqref{decor} indeed holds true, but believe these techniques to be widely applicable, which might be of independent interest. We only note that \eqref{decor} is (essentially) proven for $g(n) = {\mathfrak S}_1(n)$ and $F_x = \mathds{1}_{[0,(\log f(x))/x]}$, but we expect various other instances to work along the same lines of the proof.
In that sense, the technical aspect can be related to the above-mentioned works on Poissonian correlations of primes or sums of two squares on the real line. A main technical ingredient there is to obtain asymptotic formulae for objects such as
\[\sum_{\substack{h \leq N\\h \equiv a \mod b}} {\mathfrak S}_1(h)
\]
where $ {\mathfrak S}_1(h)$ denotes the respective singular series. The technical advancement established in this paper is to obtain asymptotics for 
\[\sum_{\substack{h \leq N\\ \lVert h\alpha\rVert \leq \log N /N}} {\mathfrak S}_1(h),
\] as long as $\alpha$ is badly approximable, and the same holds for higher-order correlations.

\subsection*{Acknowledgements}
The author would like to thank Emmanuel Kowalski and Julia Stadlmann for related discussions concerning sieve theory, and Lorenz Fr\"uhwirth for discussions on Markov chains. Additionally, the author would like to thank Nick Rome for making him aware of the work of Mikawa. Further, he would like to thank Christoph Aistleitner, Sigrid Grepstad and Niclas Technau for comments on an earlier version of this article, and the latter for making him aware of the conjecture of Larcher and Stockinger.

\section{Proof strategy and notation}\label{sec_heur}

Here we provide an overview and simplified heuristics for the proofs of the main theorems. We remark that our proof uses neither dynamical methods nor harmonic analysis and therefore clearly differs from the methods used for most previous results on Poissonian correlations mentioned in the Introduction. Instead of this, we use a combination of sieve theory, averaging techniques from analytic number theory, Ostrowski expansions from the theory of Diophantine approximation, and Markov chains.\\

For the heuristics, we  will focus on the case $k = 2$ with $R = [-1,1]$. Up to mainly combinatorial inconveniences, the cases where $k \geq 3$ as well as taking arbitrary rectangles $\bm{R}$ follows along the same lines. The result of Theorem \ref{main_thm} is first proven with a uniform roughness parameter, i.e. for the triangular arrays of rough numbers 
$(a_{n,x}^{(f)}\alpha)_{n \leq x}$. The sets of the uniform respectively non-uniform variant only differ by a negligible proportion compared to all elements in the sequence, which will be made precise in the very last step (Section \ref{sec_arr_to_seq}). 
Evaluating the pair correlations function for $(a_{n,x}^{(f)}\alpha)_{1 \leq n \leq x}$ directly, one needs to count

\[\begin{split}&\frac{1}{\Phi(x,f(x))}\#\left\{1 \leq n \neq m \leq \Phi(x,f(x)): \lVert a_n\alpha - a_m\alpha\rVert \leq \frac{1}{\Phi(x,f(x))}\right\}
= \\&\frac{1}{\Phi(x,f(x))}\sum_{\substack{1 \leq h \leq x\\ \lVert h\alpha\rVert \leq \frac{1}{\Phi(x,f(x))}}}w(h)\cdot
\#\left\{n,m \leq x: n,m \text{ $f(x)$-rough numbers}, n-m = h\right\},
\end{split}
\]
where $w(h)$ is a smooth weight function which will not cause too many problems and thus will be dropped here (see Section \ref{sec_step_fct} for how to treat this function).
Note that the expected (and in some cases provable) count differs for various $h$, depending on its arithmetical properties: As a simple showcase, observe that since $n,m$ are rough, they are in particular odd integers, and thus $h$ needs to be even in order to have any possible solutions $n - m$. A similar argument lets the count also vary for larger moduli, and one obtains for $f(x) = x^{o(1)}$ the conjectured asymptotics ($p = 2$ or in general $p \leq k$ will be treated differently)

\begin{equation}\label{sing_series_heur}\begin{split}&\#\left\{n,m \leq x: n,m \text{ $f(x)$-rough numbers}, n-m = h\right\} \sim \prod_{\substack{1 \leq p \leq f(x)\\p \nmid h}}\left(1 -\frac{2}{p}\right)\prod_{\substack{1 \leq p \leq f(x)\\p \mid h}}\left(1 -\frac{1}{p}\right)
\\&=\prod_{\substack{1 \leq p \leq f(x)}}\left(1 -\frac{2}{p}\right)\prod_{\substack{1 \leq p \leq f(x)\\p \mid h}}\left(1 -\frac{1}{p}\right)\left(1 -\frac{2}{p}\right)^{-1}.
\end{split}
\end{equation}

Indeed, by a reasonably standard procedure, using the Fundamental Lemma of Sieve Theory, one can establish this asymptotic uniformly in $h$ (see Section \ref{Sec_sieve}). We remark that one can easily adapt the method of proof used for Theorem \ref{main_thm} by deciding to not sift out certain primes, or, say, sifting out only primes with certain congruence conditions, without any serious necessary adaptations in the remainder of the proof.
\\

Since for $f(x) = x^{o(1)}$, we have
$\Phi(x,f(x)) \sim x \prod_{p \leq f(x)}\left(1 - \frac{1}{p}\right)$, and the Bohr set
\[\mathcal{B} = \mathcal{B}(\alpha,x) := \left\{1 \leq h \leq x: \lVert h\alpha \rVert \leq \frac{1}{\Phi(x,f(x))}\right\}\]
satisfies for badly approximable $\alpha$ by discrepancy arguments  (see Lemma \ref{bohr_card})
\[\#\mathcal{B} \sim 2\frac{x}{\Phi(x,f(x))} = \vol([-1,1])\frac{x}{\Phi(x,f(x))},\] we get 
\[\begin{split}&\frac{1}{\Phi(x,f(x))}\#\left\{1 \leq n \neq m \leq \Phi(x,f(x)): \lVert a_n\alpha - a_m\alpha\rVert \leq \frac{1}{\Phi(x,f(x))}\right\}
\\
\sim& \vol(R)\frac{\prod_{\substack{1 \leq p \leq f(x)}}\left(1 -\frac{2}{p}\right)}{\prod_{p \leq f(x)}\left(1 - \frac{1}{p}\right)^2}\frac{1}{\#\mathcal{B}}\sum_{h \in \mathcal{B}}\prod_{\substack{1 \leq p \leq f(x)\\p \mid h}}\left(1 -\frac{1}{p}\right)\left(1 -\frac{2}{p}\right)^{-1}.
\end{split}\]
Hence we are done if we can prove the expected asymptotic constant for $\prod_{\substack{1 \leq p \leq f(x)\\p \mid h}}\left(1 -\frac{1}{p}\right)\left(1 -\frac{2}{p}\right)^{-1}$ averaged over $\mathcal{B}$.
Ignoring for the heuristic here a convergent Euler product, this equals $g(h) := \prod_{\substack{p \leq f(x)\\p \mid h}}\left(1 + \frac{1}{p}\right)$
which is averaged over all $h \in \mathcal{B}(\alpha,x) := \{1 \leq n \leq x: \{n\alpha\} \in [0,\rho_x]\}$ for $\rho_x \approx (\log f(x))/x$. Put differently, we essentially need to establish \eqref{decor} for this $g$ and $F_x(t) := \mathds{1}_{[(\log f(x))/x,(\log f(x))/x]}(t)$.\\

To do so, we first remove the contribution of large primes since they would cause problems later. Using $g(h) \ll \frac{h}{\varphi(h)} \ll \log \log x$ and the fact that it is rare that a number $h$ is divisible by too many big primes simultaneously (see Lemma \ref{anatomy_lem}), we can replace $\prod_{\substack{p \leq f(x)\\p \mid h}}\left(1 + \frac{1}{p}\right)$ by $\prod_{\substack{p \leq \tilde{f}(x)\\p \mid h}}\left(1 + \frac{1}{p}\right)$
with $\tilde{f}$ much smaller than $f$.
After this, the average is now rewritten as \[\sum_{h \in \mathcal{B}}\sum_{\substack{d\mid h\\d\; \tilde{f}(x)\text{-smooth}}}\frac{\mu^2(d)}{d}\]
which allows us to exchange the order of summation, and we need to find asymptotics for

 \[\sum_{d}\frac{\mu^2(d)}{d}
 \frac{\#\{n \in \mathcal{B}(\alpha,x): n \equiv 0 \pmod d\}}{\#B(\alpha,x)},
 \]
where $d$ is summed over all $\tilde{f}(x)$-smooth numbers.
We aim to get the convergence to
$\sum_{d }\frac{\mu^2(d)}{d^2}$, which would lead to the desired limit value.
For constant $d \in \mathbb{N}$, the result follows from Theorem \ref{thm_bohr_div}, whereas for larger $d$, we use the structural information of Bohr sets developed mainly by Tao, Chow and Technau (see e.g.\cite{chow,CT24,tao_blog}) which allows to deduce the weaker estimate
\[\frac{\#\{n \in \mathcal{B}(\alpha,x): n \equiv 0 \pmod d\}}{\#\mathcal{B}(\alpha,x)} \ll \frac{1}{d}\] on a much larger range of $d$ (depending on $x$). This suffices to put the tail (meaning $d > T$ for a fixed large integer $T$) of the convergent series into the error term. This procedure is established in Section \ref{sec_ant_part} by applying classical tools of analytic number theory. With this at hand and leaving the proof of Theorem \ref{thm_bohr_div} for the moment aside, the proof is completed by removing the smooth weights (see Section \ref{sec_step_fct}), as well as changing from the uniform roughness of the triangular array to actual sequences (Section \ref{sec_arr_to_seq}), which follow from lengthy, but elementary computations.\\

For the proof of Theorem \ref{thm_bohr_div}, we make use of the connection between Ostrowski expansions and elements in Diophantine Bohr sets. While this is a somewhat classical topic dating back almost 100 years \cite{ostrowski}, we employ the extensive theory developed more recently in the article of Beresnevich--Haynes--Velani \cite{BHV}. For fixed $\alpha = [0;a_1,a_2\ldots,]$ written in its continued fraction expansion, we have a correspondence \[N \leftrightarrow
(b_1(N),\ldots,b_k(N)),\quad 1 \leq N < q_{k+1},\] 
where $N = \sum_{i \leq k}b_i(N)q_i, b_i \in \{0,1,\ldots,a_{i+1}\}$, and $p_k/q_k$ denotes the $k$-th convergent of $\alpha$.\\
Oversimplifying some technical arguments, we have
$N \in \mathcal{B}(\alpha,x)$ if and only if $(b_1(N),\ldots,b_k(N))$ starts with $L(x) < k(x)$ many $0$'s, and the assumption of $x\rho_x \to \infty$ ensures $k(x) - L(x) \to \infty$ (In fact, this would only work for $x = q_{k+1}$ and $\rho_x = \lVert q_{L}\alpha\rVert$. Therefore, we pass to a subsequence $(x_n)_{n \in \mathbb{N}}$ where $x_{n+1}/x_n \to 1$ with $x_n = \sum_{i = 0}^{M}b_{i}q_{k_n-i}$ for arbitrarily slowly growing $M(x)$, and approximate the corresponding $\rho_x$ by $\sum_{i = 0}^{M}b_{i}\lVert q_{L+i}\alpha\rVert$. As long as we can ensure $k(x) - L(x) - 2M(x) \to \infty$, we are essentially back in our simplified situation.) If we now view $\{1,\ldots,q_{k+1}-1\}$ as a probability space with uniform probability measure, each $b_i \pmod d, L(x) < i \leq k(x)$ can be viewed as a random variable on $\{0,\ldots,d-1\}$. If these $b_i$ had a reasonably non-degenerated distribution with $b_i$ being stochastically independent of the other $b_j$'s, the question would almost reduce to asking for the probability that after $K$ dice rolls, the sum of the dots shown being $a \pmod d$, when $a,d$ are fixed and $K \to \infty$. Clearly, one expects an equidistribution for this question, and this can be proven by a standard application of the theory of Markov chains.\\

However, the random variables $b_i$ are not independent, but follow the Markov chain condition $b_{i} = a_{i+1} \Rightarrow b_{i-1} = 0$ (we remark that there are two different Markov chains here: One for the Ostrowski digits, and one coming from the random walk modulo $d$). The idea to circumvent this obstruction is to observe that conditioned on the event 
$b_{i+3} = b_{i} = 0$, the random variable $(b_{i+1},b_{i+2})$ is now independent of all the remaining $b_j, j \notin \{i,i+1,i+2,i+3\}$. The reason we take pairs is to exploit the fact that successive convergent denominators are always coprime, and thus summing the independent random variables $b_iq_i + b_{i+1}q_{i+1} \pmod d$ gives rise to a Markov chain that is irreducible (picking just an arbitrary collection of $i$'s might e.g. have resulted in $q_i$ always being even, destroying any hope of establishing equidistribution $\pmod d$ for even $d$).
It turns out that the occurrence of such a constellation $(0,*,*,0)$ for an unbounded number of $i$'s is generic among all tuples: This is essentially just arguing that flipping a coin infinitely often, we get almost surely infinitely often ``heads'', since the event $b_i = 0$ can be shown to have a positive probability, regardless of any given values for $b_j, j \neq i$. In that way, we get our desired independence assumption, and can thus appeal to the Markov chain convergence described above. This whole procedure is executed in Section \ref{sec_mod_d}.\\
We remark that there might be a possibility to prove Theorem \ref{thm_bohr_div} by doing a lattice-point counting argument in a similar way as in \cite{TW20} building upon the work of Skriganov \cite{Skri}, as remarked by N. Technau in private communication. However, we decided to opt for the approach using Markov chains sketched above, which provides a new perspective on the matter.\\

While the property of $\alpha$ being badly approximable is used on numerous occasions in the proof, the most important reason why Theorem \ref{main_thm} fails for well-approximable numbers (or in other words, why Theorem \ref{dioph_ass_thm} holds) is that in this case, the cardinality of $\mathcal{B}(\alpha,x)$ does \textit{not} follow the expected asymptotic $\#\mathcal{B}(\alpha,x) \sim \vol(R) \frac{x}{\Phi(x,f(x))}$:
If $p_n/q_n$ approximates $\alpha$ very well, then choosing $x$ a suitable multiple of $q_n$ provides too many $1 \leq h \leq x$ to lie in $\mathcal{B}(\alpha,x)$.

Note that \[\prod_{\substack{2 < p \leq f(x)\\p \mid h}}\left(1 -\frac{1}{p}\right)\left(1 -\frac{2}{p}\right)^{-1} \gg 1,\] which is of the same order as the expected value, making it unnecessary to execute the delicate averaging technique as applied for Theorem \ref{main_thm}. 
Thus the desired limit is off by a factor of order
$\frac{\#\mathcal{B}(\alpha,x)}{\vol(\bm{R}) \frac{x}{\Phi(x,f(x))}},$ which along a subsequence is unbounded.
This is essentially the proof of Theorem \ref{dioph_ass_thm}.

\subsubsection*{Notation}
We use the Vinogradov notation $\ll$ as well as the $O$- and $o$-notation. We use $f\sim g$ for $\lim_{n \to \infty} \frac{f(n)}{g(n)} = 1$ and $f \asymp g$ for $f \ll g$ and $g \ll f$. We use the standard notations $\mu,\tau,\omega,\varphi$ for Möbius function, number of divisors function, number of (distinct) prime divisors function, and Euler's totient function, respectively. When the summation respectively product is over $p$, we interpret this as a sum respectively a product over prime numbers.
We denote by $\{x\}$ the integer part of $x$ and by $\lVert x\rVert$ the distance to the nearest integer.
\\
We introduce some vector notation that will be used throughout this article: We write vectors $\bm{n} \in \mathbb{Z}^k$ in bold and integers $n \in \mathbb{Z}$ in normal font.
For a product set $\bm{S} = A_1 \times A_2 \times \ldots \times A_k$, we denote by
$\bm{S}_{\neq}$ the set of all tuples $(a_1,a_2,\ldots,a_k) \in \bm{S}$ with pairwise distinct coordinates.
For an integer vector $\bm{n} = (n_1,\ldots,n_k) \in \Z^k$, we denote by $n^{+} := \max_{1 \leq i \leq k}n_i,$ $n^{-} := \min_{1 \leq i \leq k}n_i$.
Further, we denote by $P^{-}(\bm{n}) := \min_{1 \leq i \leq k} P^{-}({n_i})$ where $P^{-}(n)$ denotes the smallest prime dividing $n$. 
Moreover, for $t \in \mathbb{\bm{R}}$, we let $t_{>0} := \max\{t,0\}$ and $t_{<0} := \min\{t,0\}$.

\section{Sieve estimate}\label{Sec_sieve}

In this section, we provide the sieve-theoretic application that allows to establish the analogue of the Hardy-Littlewood $k$-admissibility conjecture for rough numbers, provided that the roughness parameter is $x^{o(1)}$:

\begin{lem}\label{sieve_k_tuple}
Let $f(x) = x^{r(x)}$ where $r(x)$ decreases monotonically to $0$.
Let $k \geq 2$ be a fixed integer and for $\bm{h} \in \Z^{k-1}_{\neq}$, let
$g_{\bm{h}}(p) = \#\{0 \pmod p,h_{1} \pmod p,\ldots,h_{k-1} \pmod p\}$.
Then we have
\[\begin{split}&\#\{(m,\bm{n}) \in [1,x]^{k}_{\neq}: P^{-}((m,\bm{n})) > f(x), m\bm{1}_{k-1} - \bm{n} = \bm{h}\} \\\sim&\;
(x - (h^{+}_{>0}-h^{-}_{<0}))_{>0}\prod_{1 \leq p \leq f(x)}\left(1 - \frac{g_{\bm{h}}(p)}{p}\right) + O(x^{o(1)}),\end{split}\]
with the convergence being uniform in $\bm{h}$.
\end{lem}

For the proof, we will make use of the Fundamental Lemma of Sieve Theory in the form of \cite[Theorem 18.11]{kouk_book} which we state below for the reader's convenience.

\begin{lem}\label{fund_lem}
    Let $(b_n)_{n \geq 1}$ be non-negative reals such that 
    $\sum_{n=1}^{\infty}b_n < \infty.$ 
 Let $\mathcal{P}$ be a finite set of primes, and write 
 $P = \prod_{p \in \mathcal{P}} p$. Set $y = \max \mathcal{P}$,
and 
$B_d = \sum_{n \equiv 0 \mod d}b_n$.
Assume that there exists a multiplicative function $g$ such that
$g(p) < p$ for all $p \in \mathcal{P}$, a real number $X$, and positive constants $\kappa,C$ such that

\begin{equation}\label{A_d}B_d =: X\frac{g(d)}{d} + r_d, \quad d | P,\end{equation}
and

\begin{equation}\label{sieve_dim}\prod_{p \in (y_1,y_2] \cap \mathcal{P}}\left(1 - \frac{g(p)}{p}\right) < \left(\frac{\log y_2}{\log y_1}\right)^{\kappa}\left(1 + \frac{C}{\log y_1}\right), \quad 3/2 \leq y_1 \leq y_2 \leq y.\end{equation}

Then, uniformly in $u \geq 1$ we have
\[\sum_{(n,P) = 1}b_n = \left(1 + O_{\kappa,C}(u^{-u/2})\right)X \prod_{p \in \mathcal{P}}\left(1 - \frac{g(p)}{p}\right) + O\left(\sum_{\substack{d \leq y^u\\d \mid \mathcal{P}}} |r_d|\right).\]
\end{lem}

\begin{proof}[Proof of Lemma \ref{sieve_k_tuple}]
First, we translate the above count into the setup of Lemma \ref{fund_lem}. We observe that in order to satisfy $m\bm{1}_{k-1} - \bm{n} = \bm{h}$ for some $(m,\bm{n}) \in [1,x]^{k}_{\neq}$, we have
\[\max\{0,-h^{-}\} \leq m \leq x - \max\{h^{+},0\}.\]
Thus we define \[b_n = b_n^{(\bm{h})} := \sum_{-h^{-}_{>0} \leq m \leq x - h^{+}_{>0}}\mathds{1}_{[m\prod_{i = 1}^{k-1}(m+h_i) = n]},\]
and observe that 
\[
\#\{(m,\bm{n}) \in [1,x]^{k}_{\neq}: P^{-}((m,\bm{n})) > f(x), m\bm{1}_{k-1} - \bm{n} = \bm{h}\}
= \sum_{\substack{n \in \N\\(n,P) = 1}} 
b_n
\]
where $P := \prod_{p \leq f(x)}p$.
    This approach itself is not novel, see e.g. \cite[Example 18.2]{kouk_book} or \cite[Theorem 6.7]{iwaniec_kowalski} for related setups.\\

We may assume that for $p \leq k$ we have $g_{\bm{h}}(p) < p$, since otherwise there are no solutions at all and the result follows trivially.
In order to show \eqref{A_d}, observe that on full periods $\mathbb{Z}_d$ in the sieve range $-(h^{-}_{<0}) \leq m \leq x - h^{+}_{>0}$, we get precisely $g_{\bm{h}}(d) := \prod_{p \mid d}g_{\bm{h}}(p)$ many solutions, which proves \eqref{A_d}
with $X = (x - (h^{+}_{>0}-h^{-}_{<0}))_{>0}, g = g_{\bm{h}}$, and $|r_d| \leq g_{\bm{h}}(d) \leq k^{\omega(d)}$.

Thus an application of Lemma \ref{fund_lem} with $\kappa = k, u = \frac{1}{\sqrt{r(x)}}, y = f(x)$
where $f(x) = x^{r(x)}$, and the well-known estimate $k^{\omega(d)} \ll_{\varepsilon,k} d^{\varepsilon}$ shows that

\[\sum_{\substack{n \in \N\\(n,P) = 1}} b_n
\sim (x - (h^{+}_{>0}-h^{-}_{<0}))_{>0} \prod_{1 < p \leq f(x)}\left(1 - \frac{g_{\bm{h}}(p)}{p}\right) + O\left(f(x)^{\frac{1}{\sqrt{r(x)}}+1}\right),
\]
with the convergence being uniform in $\bm{h}$.
\end{proof}

\section{Equidistribution mod $d$ in Bohr sets}\label{sec_mod_d}

\subsubsection{Prerequisites - Ostrowski expansions}
For the convenience of the reader, we recall here some well-known facts on Diophantine approximation and continued fractions that are used in this paper. For a more detailed background, see the classical literature e.g. \cite{all_shall,rock_sz,schmidt}. Furthermore, we make use of various results of \cite{BHV} where Bohr sets in Diophantine approximation were developed, which we will recall and slightly generalize to the signed setup.\\

Every irrational $\alpha$ has a unique infinite continued fraction expansion $[a_0;a_1,...]$ with  convergents $p_k/q_k = [a_0;a_1,...,a_k]$ satisfying the recursive formulae
\begin{equation*}
p_{k+1} = p_{k+1}(\alpha) = a_{k+1}p_k + p_{k-1}, \quad q_{k+1} = q_{k+1}(\alpha) = a_{k+1}q_k + q_{k-1}, \quad k \gs 1\end{equation*}
with initial values $p_0 = a_0,\; p_1 = a_1a_0 +1,\; q_0 = 1,\; q_1 = a_1$.
One can deduce from these recursions that $q_k$ grows exponentially in $k$; in particular,
we have for any $k,j \in \mathbb{N}$
\begin{equation}\label{exp_grow}
    \frac{q_{k+j}}{q_{k}} \gg (3/2)^j,
\end{equation}
with an absolute implied constant.
Defining 
\[\delta_k := q_k \alpha - p_k,\quad k \geq 0,\]
we have the following identities for any $k \geq 1:$

\begin{align}
    \label{id1}\delta_0 &= \{\alpha\}, \quad \delta_k = (-1)^k\lVert q_k\alpha\rVert,\\
    a_{k+1}\delta_k &= \delta_{k+1} - \delta_{k-1},\\
    \frac{1}{2} &\leq q_{k+1}|\delta_k|\leq 1,\\
    \left\lvert\frac{\delta_{m+k}}{\delta_m}\right\rvert& \ll \left(\frac{2}{3}\right)^k,\label{exp_decay}\\
    |\delta_k| &= a_{k+2}|\delta_{k+1}| + |\delta_{k+2}| ,\label{id_23rd_last}\\
    |\delta_k| & = \sum_{i = 1}^{\infty}a_{k+2i}|\delta_{k+2i-1}|,\label{id_2nd_last}\\
    |\delta_k| + |\delta_{k+1}| &= \sum_{i = k+1}^{\infty}a_{i+1}|\delta_i|.\label{id_last}
\end{align}

We call an irrational number $\alpha$ \textit{badly approximable} if $C(\alpha) := \sup_{n \in \N} a_n(\alpha) < \infty$. 
Fixing an irrational $\alpha = [a_0;a_1,...]$, the Ostrowski expansion of a non-negative integer $N$ is the unique representation

\begin{equation*}
N = \sum_{i = 0}^k b_{i}(N)q_{i}, \quad
b_{k+1} \neq 0,\;\; 0 \ls b_0 < a_1, \;\;  0 \ls b_{i} \ls a_{i+1} \text { for } i \gs 1,
\end{equation*}
with the additional rule that
$b_{i-1} = 0$ whenever $b_{i} = a_{i+1}$.
If $b_0(N) = b_1(N) = 0$, we obtain the very useful identity
\begin{equation}\label{Ostr_id}\{N\alpha\} = \sum_{i = 0}^k b_{i}(N)\delta_{i},\end{equation}
from which we can deduce that 
\[(b_{m}-1)|\delta_m| + (a_{m+2} - b_{m+1})|\delta_{m+1}| \leq \lVert N\alpha\rVert \leq (b_{m}+1)|\delta_m|,\]
where $m$ is the smallest index such that $b_m(N) \neq 0$, see \cite[Lemma 4.1 and Lemma 4.2]{BHV}.
\\

We now recall the classical Denjoy--Koksma inequality from Discrepancy theory (see e.g. \cite[Chapter 2, Theorem 3.4]{KN}):

\begin{prop}[Denjoy-Koksma inequality]\label{denj_koks}
    Let $f: \mathbb{R} \to \mathbb{R}$ be a $1$-periodic function with bounded total variation $Var(f)$ on $[0,1)$ and let $N = \sum_{i = 0}^k b_iq_i$ be an integer written in Ostrowski expansion with respect to some irrational $\alpha$. Then 
    \[\sum_{n \leq N}f(n\alpha) = N \int_{0}^1 f(x)\,\mathrm{d}x + O\left(Var(f) \sum_{i = 0}^k b_i\right),\]
    with an absolute implied constant.
\end{prop}

We say that a tuple 
$(c_{i_0},\ldots,c_{i_r}), i_t \in \N$ is \textit{admissible} (with respect to $\alpha$) if
there exists $N \in \N$ such that\footnote{Note the abuse of notation: Formally, we would need to define $((c_{i_0},i_0),\ldots,(c_{i_r},i_r))$ to be admissible but this rather clumsy notation is omitted for the sake of readability.}
\begin{equation}\label{admissible_def}N = \sum_{i = 0}^k b_{i}q_{i},\quad  b_{i_t} = c_{i_t}, 1 \leq t \leq r,\end{equation}
where $N$ is written in its Ostrowski expansion.
Given an admissible tuple $(d_0,\ldots,d_m)$, we define the \textit{Ostrowski cylinder sets} as
    
    \[A(d_0,\ldots,d_m) := \{N \in \N: b_i(N) = d_i\quad \forall 1 \leq i \leq m\}, \quad A_M(d_0,\ldots,d_m) := A(d_0,\ldots,d_m) \cap [1,M].\] 
    We say $(c_{i_0},\ldots,c_{i_r}), i_t \in \N$ is \textit{admissible} in 
$A(d_0,\ldots,d_m)$ respectively in $A_M(d_0,\ldots,d_m)$ if there exists $N \in A(d_0,\ldots,d_m)$
respectively $N \in A_M(d_0,\ldots,d_m)$ such that \eqref{admissible_def} holds.\\

We define the Ostrowski expansion for real numbers as follows: Fixing $\alpha$, the \textit{real} Ostrowski expansion of a real number $\gamma \in [-\alpha,1 - \alpha)$ is the unique representation 

\begin{equation}\label{real_ostr}\gamma = \sum_{i = 0}^{\infty} c_{i}(\gamma)\delta_i \quad \text{ where }
 0 \ls b_0 < a_1, \;\;  0 \ls b_{i} \ls a_{i+1} \text { for } i \gs 1,\end{equation}
with the additional rule that
$c_{i-1} = 0$ whenever $c_{i} = a_{i+1}$.

\begin{prop}\label{prop_N_pm}
    Let $\alpha$ be an irrational number and let
    \[x = \sum_{i \leq n}b_i(x)q_i, \quad b_n\neq 0\]
    be an integer written in its Ostrowski expansion.
    For $m \in \mathbb{N}$ with $m < n$, let us denote 
    $x_m^{-} := \sum_{n-m \leq i \leq n}b_i(x)q_i$
    and let 
    $x_m^{+} = \sum_{i \geq n-m}b_i(x_m)q_i$
    where $x_m = x + q_m$. Then we have the following:
    \begin{itemize}
    \item[(i)] The respective Ostrowski expansions of both $x_{m}^{-}, x_{m}^{+}$ have at most $m+2$ non-zero digits. 
        \item[(ii)] $x_{m}^{-} \leq x \leq x_m^{+}$.
        \item[(iii)] $\frac{x_m^{+}}{x_m^{-}} \leq 1 + O\left(\left(\frac{2}{3}\right)^m\right)$,
        with an absolute implied constant.
    \end{itemize}
\end{prop}

\begin{proof}
    Since $x_{m}^{-}, x_{m}^{+}$ arise from truncations of Ostrowski expansions of $x$ and $x_m$ respectively, and 
    $x < x_m < q_{n+2}$, (i) follows immediately. While $x_{m}^{-} \leq x$ is trivial, note that 
    \begin{equation}\label{tail_ostr}
    x_m^{+} = x + q_m - \underbrace{\sum_{i < m}b_i(x_m)q_i}_{< q_m} > x,\end{equation}
    thus (ii) follows. By the same argumentation, 
    we obtain
    $x_m^{-} > x - q_m$, and since $x > q_n$, 
    \[1 \leq \frac{x_m^{+}}{x_m^{-}} \leq \frac{x + q_m}{x - q_m} = 1 + O(q_{n-m}/q_n)
    = 1 + O\left(\left(\tfrac{2}{3}\right)^m\right),\]
    using \eqref{exp_grow}.
\end{proof}

We also need the ``dual'' of Proposition \ref{prop_N_pm} that considers upper and lower bounds for truncated real Ostrowski expansions.
\begin{prop}\label{dual_real}
    Let $\gamma \in [-\alpha,1-\alpha)$ with 
    $\gamma = \sum_{i = m}^{\infty} c_{i}(\gamma)\delta_i, c_m(\gamma) \neq 0$  being its real Ostrowski expansion and assume $m \geq 1$. Further, let $L \in \mathbb{N}$ be given and define 
    \[\gamma^+_{L,m} = \sum_{i \leq L+m}c_{i}(\gamma^+)\delta_i, \quad \gamma^-_{L,m} = \sum_{i \leq L+m}c_{i}(\gamma^-)\delta_i\]
    where $\gamma^+ = \gamma + |\delta_{m+L}|, \gamma^- = \gamma - |\delta_{m+L}|$. Then we have the following:
    \begin{itemize}
        \item[(i)] $\gamma^-_{L,m} \leq \gamma \leq \gamma^+_{L,m}$.
        \item[(ii)] $\frac{\gamma^+_{L,m}}{\gamma^-_{L,m}} \leq 1 + O((\tfrac{2}{3})^L)$,
        with an absolute implied constant.
          \item[(iii)] The Ostrowski expansions of $\gamma^+_{L,m}, \gamma_{L,m}^-$ are of the form $(\underbrace{0,\ldots,0}_{m -2},d_{m-1},\ldots,d_{m+L})$.
        \end{itemize}
\end{prop}

\begin{proof}
From \eqref{id_2nd_last} and \eqref{id_last} we infer that 
\begin{equation}\label{first_ostr_det}|\delta_{m+1}| \leq |\gamma| \leq |\delta_{m-1}|.\end{equation}
    From \eqref{id_2nd_last} we obtain that for any $\beta = \sum_{i = 1}^{\infty} \tilde{c}_{i}\delta_i$ and any $k \geq 1$ that
    $\left|\sum_{i > k} \tilde{c}_{i}\delta_i\right| \leq |\delta_{k}|.$
    Therefore, we obtain \[\gamma_{L,m}^+ \geq \gamma + |\delta_{m+L}| - \left\lvert\sum_{{i > L+m}}c_{i}(\gamma^+)\delta_i\right\rvert \geq \gamma,
    \]
    and analogously, $\gamma^-_{L,m} \leq \gamma$, showing (ii). By the same argumentation, we have
    \[\gamma_{L,m}^+ \leq \gamma + |\delta_{m+L}| + \left\lvert\sum_{{i > L+m}}c_{i}\right\rvert \leq 2|\delta_{m+L}|\]
    as well as $\gamma_{L,m}^- \geq \gamma - 2|\delta_{m+L}|$ which in view of 
    \eqref{exp_decay} and \eqref{first_ostr_det} proves (ii). Since 
    $|\gamma_{L,m}^{\pm}| \leq|\delta_{m+L}| + |\delta_{m-1}|$, using 
    again \eqref{first_ostr_det} shows (iii).
\end{proof}

\subsection{Prerequisites - Markov chains and conditional probabilities}

We will make use of the very basic theory of Markov chains in order to obtain the equidistribution in $\mathbb{Z}_d$. For more details we refer the interested reader to the standard literature on this topic such as \cite{MeyTwe,Woess}.\\

\begin{defi}
Given a probability space $(\Omega,\mathcal{A},\mathbb{P})$, a (time-homogeneous) Markov chain with finite state space $X$ is a sequence of random variables $(Z_i)_{i \in \mathbb{N}_0}$ with $Z_i: \Omega \to X$ that satisfies the following properties:

\begin{itemize}
    \item[(i)] For all $x_0,\ldots,x_{n+1} \in X, n \in \mathbb{N}_0$ such that 
    $\mathbb{P}[Z_n = x_n,\ldots, Z_0 = x_0] > 0$, one has
    \[\mathbb{P}[Z_{n+1}= x_{n+1} \mid Z_n = x_n,\ldots, Z_0 = x_0] = \mathbb{P}[Z_{n+1}= x_{n+1} \mid Z_n = x_n].\]
    \item[(ii)] For all $x,y \in X$ and all $n,m \in \mathbb{N}_0$ such that
    $\mathbb{P}[Z_n = x]> 0, \mathbb{P}[Z_m = x]> 0$, we have
    \[\mathbb{P}[Z_{n+1} = y \mid Z_n = x] = \mathbb{P}[Z_{m+1} = y \mid Z_m = x].\]
\end{itemize}
We define the transition matrix $P$ by $P(x,y) := \mathbb{P}[Z_{n+1} = y \mid Z_n = x]$ and write 
$P^{(k)}(x,y) = \mathbb{P}[Z_{n+k} = y \mid Z_n = x]$. The initial distribution of  $(Z_i)_{i \in \mathbb{N}_0}$ is the distribution $\nu$ on $X$ defined by
$\nu(x) := \mathbb{P}[Z_0 = x], x \in X$.\\

We call a Markov chain \textbf{irreducible} when for all $x,y \in X$, there exists $n \in \mathbb{N}$ such that $P^{(n)}(x,y) > 0$. We call a Markov chain \textbf{aperiodic} if for all $x \in X$
\[\gcd(\{n > 0: P^n(x,x) > 0\}) = 1.\]
We call a distribution $\pi$ on $X$ \textbf{stationary} if (interpreted as vector $\pi_x := \pi(x), x \in X$) we have $\pi P = \pi$.
\end{defi}

With these definitions at hand, we can state the Markov chain convergence Theorem with finite state space for time-homogeneous chains in the form of \cite[Theorem 3.28]{Woess}:

\begin{prop}\label{prop_woess}
    Let $(Z_i)_{i \in \mathbb{N}_0}$ be an irreducible, aperiodic time-homogeneous Markov chain with finite state space $X$. Then there are $k \in \mathbb{N}$ and $\tau < 1$ such that for the unique stationary distribution $\pi$ one has
    \begin{equation}
        \sum_{y \in X}\left|P^{(n)}(x,y) - \pi(y)\right| \leq 2\tau^n
    \end{equation}
    for every $x \in X$ and $n \geq k$.
\end{prop}

We will apply this proposition to the random walk on the additive groups $\mathbb{Z}_d$ in the following way:

\begin{cor}\label{rand_walk}
Let $d \in \mathbb{N}$, and $\mu$ a probability measure on $\mathbb{Z}_d = \{0,\ldots,d-1\}$ that satisfies:
\begin{itemize}
    \item[(i)] $\mu(0) > 0$.
    \item[(ii)] There exists $a,b \in \mathbb{Z}_d$ with  $\mu(a) > 0, \mu(b) > 0$ such that $\gcd(a,b,d) = 1$.
\end{itemize}
Given a probability space $(\Omega,\mathcal{A},\mathbb{P})$, let $(Y_i)_{i \in \mathbb{N}_0}$ be a sequence of independent random variables on $\mathbb{Z}_d$, with $Y_0 \stackrel{d}{=} \theta$ and $Y_i \stackrel{d}{=} \mu, i \geq 1$. Then we have the following:
For all $\varepsilon > 0$, there exists $K_0 = K_0(\mu,\varepsilon) \in \mathbb{N}$ such that for all distributions $\theta$ on $\mathbb{Z}_d$, all $a \in \mathbb{Z}_d$ and all $n \geq K_0$, we have
\[\left\lvert\mathbb{P}\left[\sum_{i = 0}^n Y_i = a\right] - \frac{1}{d}\right\rvert < \varepsilon.\]
\end{cor}

\begin{proof}
A routine computation shows that the partial sums $Z_n := \sum_{0 \leq i \leq n}Y_i$ form a time-homogeneous Markov chain with initial distribution $Y_0$ and transition probabilities
$P(x,y) := \mu(y-x)$.
By employing Bezout's identity on $a,b,d$, we find under (ii) the existence of $k,\ell \in \mathbb{N}$ such that $ak + b\ell \equiv 1 \pmod d$. Thus $P^{(k+\ell)}(0,1) > 0$, and hence the Markov chain is irreducible. Further, (i) immediately implies the aperiodicity. Note that the transition matrix $P$ with 
    $P(i,j) := \mu(j-i), i,j \in \mathbb{Z}_d$ is doubly stochastic, and thus the unique stationary measure is $\pi = (\tfrac{1}{d},\ldots,\tfrac{1}{d})$. An application of Proposition \ref{prop_woess} now proves the claim.
\end{proof}

Finally, we recall some standard notation for conditional probabilities that will be used in this section. We remark that we will restrict ourselves to random variables with finite state space.
Given a probability space $(\Omega,\mathcal{A},\mathbb{P})$, and a random variable $Z: \Omega \to S$
to some (finite) space $S$, we define for $A \in\mathcal{A}$ with $\mathbb{P}[A] > 0$ the random variable
$Y = Z| A$ via \[\mathbb{P}[Y = a] := \mathbb{P}[Z = a|A] = \frac{\mathbb{P}[(Z = a) \cap A]}{\mathbb{P}[A]}, \quad a \in S.\]
For a Sigma-algebra $\mathcal{B} \subseteq \mathcal{A}$, we denote by
$\mathbb{E}[Y|\mathcal{B}]$ the conditional expectation, and for an event $A \in \mathcal{A}$, we write $\mathbb{P}[A|\mathcal{B}] := \mathbb{E}[\mathds{1}_{A}|\mathcal{B}]$.
    As usual, for random variables $(X_1,\ldots,X_n)$ we write
    $\mathbb{E}[Y|(X_1,\ldots,X_n)] := \mathbb{E}[Y|\sigma((X_1,\ldots,X_n))]$
    where $\sigma(X)$ denotes the smallest Sigma-algebra generated by $X$. 

\subsection{Proof of Theorem \ref{thm_bohr_div}}

\begin{lem}[Bohr set cardinality]\label{bohr_card}
Let $(I_x)_{x \in \mathbb{N}}$ a sequence of intervals in $[0,1]$ with $\lim_{x \to \infty }x\lambda(I_x) = \infty$, let $\alpha$ be badly approximable and
define
    \[\mathcal{B}(x,I_x) = \mathcal{B}_{\alpha}(x,I_x) := \{1 \leq n \leq x: \{n\alpha\} \in I_x\}.\]
Then we get
\begin{equation}\label{bohr_asymp}\#\mathcal{B}(x,I_x) \sim \lambda(I_x)I_x, \quad x \to \infty.\end{equation}
\end{lem}

\begin{proof}
Let $x^-_m \leq x \leq x^+_m$ be defined as in Proposition \ref{prop_N_pm} which we now apply.
    By (ii), we obtain
    \[\#\mathcal{B}(x_{m}^{-},I_x)\leq \#\mathcal{B}(x,I_x) \leq \#\mathcal{B}(x_{m}^{+},I_x).\]
    By (iii), we observe that $\sum_{i}b_i(x_m^+) \ll_{\alpha} \frac{m}{x_m^+}$, thus 
    by the Denjoy--Koksma inequality (Proposition \ref{denj_koks}), we get
    \[\#\mathcal{B}(x_{m}^{+},I_x) = \sum_{n \leq x_m^{+}} \mathds{1}_{[\{n\alpha\} \in I_x]}
    = \lambda(I_x)x_{m}^{+} + O_{\alpha}(m),\]
    and the same argument shows
    \[\#\mathcal{B}(x_{m}^{-},I_x) = \lambda(I_x)x_{m}^{-} + O_{\alpha}(m).\]
    Setting $m = \lfloor\sqrt{x\lambda(I_x)}\rfloor$ and using $q_{j-m}/q_j \ll (2/3)^m$, \eqref{bohr_asymp} follows from (ii).
\end{proof}

\begin{prop}\label{HK_machine}
Let $\mathcal{B}_{\alpha}(x,\rho) = 
\{1 \leq n \leq x: \lVert n\alpha\rVert < \rho\}$ 
and let $d \in \mathbb{N}$ and $j \in \mathbb{Z}_d$ arbitrary.
Then for any badly approximable $\alpha$, we have
    \[\sum_{\substack{h \in \mathcal{B}(x,\rho)\\h \equiv j \pmod d}}1
    \ll_{\alpha} \frac{\#\mathcal{B}_{\alpha}(x,\rho)}{d} + \sqrt{\#\mathcal{B}_{\alpha}(x,\rho)},
    \]
    with the implied constant depending only on $\alpha$.
\end{prop}

\begin{proof}
    This follows by obvious modifications from the proof of \cite[Lemma 7.4]{HK}, 
    using crucially the fact that \[\mathcal{B}_{\alpha}(x,\rho) \subseteq
    \{au + bv: \lvert a \rvert \leq z, \lvert b \rvert \leq z\}\]
    for some $z$ satisfying $z^2 \ll_{\alpha} \#\mathcal{B}_{\alpha}(x,\rho)$
    (see e.g. the combination of \cite[Lemmas 2.7 and Lemma 2.15]{CT24}, or \cite{tao_blog} for a longer discussion on these Bohr set inclusions).
\end{proof}

While Proposition \ref{HK_machine} is a powerful tool for error terms, the implied constant is unavoidable, and thus there is no way that the employed method could provide an asymptotics as stated in Theorem \ref{thm_bohr_div}. To this end, we need the following refinements of Proposition \ref{HK_machine} in the two subsequent statements that become slightly technical.

\begin{prop}[Approximating Bohr sets by Ostrowski cylinder sets]\label{cylinder_cover}
    Let $\alpha$ be badly approximable with $C(\alpha) = \sup_{i \in \mathbb{N}}a_i < \infty$ and $I = [0,\gamma]$ or $I = [1 - \gamma,1], \gamma > 0$ be an interval with $|\delta_m| \leq \gamma < |\delta_{m+1}|$, $m \geq 2$.
    Let $x,L \in \mathbb{N}$ such that $x\gamma > (C(\alpha)+1)^{L}$.Then we have the following:

There exists sets $S^+,S^{-}$ with $|S^+|, |S^-| \ll L^{C(\alpha)+1}$ such that 
    \begin{equation}\label{bohr_sandwich}\bigcup_{(\underbrace{0,\ldots, 0}_{m-1}, d_{m-1},\ldots,d_{m+L}) \in S^-} A_x(0,\ldots,d_{m+L}) \subseteq B_{\alpha}(x,I) \subseteq
\bigcup_{(\underbrace{0,\ldots, 0}_{m-1}, d_{m-1},\ldots,d_{m+L}) \in S^+} A_x(0,\ldots,d_{m+L}),
\end{equation}
with 
\begin{equation}\label{sandwich_error}\Bigg|\bigcup_{(\underbrace{0,\ldots, 0}_{m-1}, d_{m-1},\ldots,d_{m+L}) \in S^+} A_x(0,\ldots,d_{m+L}) \setminus \bigcup_{(\underbrace{0,\ldots, 0}_{m-1}, d_{m-1},\ldots,d_{m+L}) \in S^-} A_x(0, \ldots,d_{m+L})\Bigg| \ll_\alpha {x\gamma\left(\tfrac{3}{2}\right)^{-L}}.\end{equation}
\end{prop}

\begin{proof}
    We will only prove the case where $I = [0,\gamma]$; the case $I = [1 - \gamma,1]$ works in the completely same fashion by alternating signs adequately. Furthermore, we may assume that $\rho \in [-\alpha,1-\alpha)$ since we otherwise simply add a suitable integer. We now define the following partial order on integer-valued tuples respectively sequences: 
    Let $(b_0,\ldots,b_i) \in \mathbb{Z}^i\setminus \{(0,\ldots,0)\}, i \in \mathbb{N} \cup \{\infty\}$ and let $i_0 = \min\{0 \leq t \leq i: b_t \neq 0\}$. 
    Then 
    \[(0,\ldots,0) \preceq (b_0,\ldots,b_i) :\Leftrightarrow (-1)^{i_0}b_{i_0} > 0.\]
        Letting
$(b_0,\ldots,b_i),(b_0',\ldots,b_j')$ be such tuples we define
    \[(b_0,\ldots,b_i) \preceq (b_0',\ldots,b_j') :\Leftrightarrow (0,\ldots,0) \preceq (b_0'-b_0,\ldots,b'_{\min\{i,j\}}-b_{\min\{i,j\}}).
    \]
    Using
    \eqref{Ostr_id} and the identities \eqref{id1} - \eqref{id_last}, it is immediate that
    \begin{equation}\label{bohr_equiv}N \in \mathcal{B}_{\alpha}(x,I) \Leftrightarrow
    (0,0,\ldots) \preceq (b_0(N),b_1(N),\ldots) \preceq (c_0(\gamma),c_1(\gamma),\ldots)
    \end{equation}
    where 
    $N = \sum_{i = 0}^{\infty} b_{i}(N)q_{i}, \quad \gamma = \sum_{i = 0}^{\infty} c_{i}(\gamma)\delta_i$
    are both written in their respective Ostrowski expansions.\\

We now define $\gamma_{L,m}^+,\gamma_{L,m}^-$ as in Proposition \ref{dual_real}
and let
\[\begin{split}S^+ &:= \left\{(d_0,\ldots,d_{m+L}) \text { admissible}: (0,0,\ldots) \preceq (d_0,\ldots,d_{m+L})
\preceq (c_0(\gamma^+_{m+L}),\ldots,c_{m+L}(\gamma^+_{m+L}))
\right\},\\
S^- &:= \left\{(d_0,\ldots,d_{m+L}) \text { admissible}: (0,0,\ldots) \preceq (d_0,\ldots,d_{m+L})
\preceq (c_0(\gamma_{m+L}^{-}),\ldots,c_{m+L}(\gamma_{m+L}^{-}))
\right\}.\end{split}
\]
By \eqref{bohr_equiv} and Proposition \ref{dual_real} (i), we see that

\[\bigcup_{(d_0,\ldots,d_{m+L}) \in S^-} A_x(d_0,\ldots,d_{m+L}) \subseteq \mathcal{B}_{\alpha}(x,I) \subseteq
\bigcup_{(d_0,\ldots,d_{m+L}) \in S^+} A_x(d_0,\ldots,d_{m+L}).
\]
By Proposition \ref{dual_real} (iii),  we have $(d_0,\ldots,d_{m+L}) = (0,\ldots,0,d_{m-1},\ldots,d_{m+L})$ for all $(d_0,\ldots,d_{m+L}) \in S^- \cup S^+$. Since $d_i \in \{0,\ldots,a_{i+1}\} \subseteq \{0,\ldots,C(\alpha)\}$, we conclude that $|S^+|, |S^-| \ll L^{C(\alpha)+1}$.\\
Furthermore, observe that by Proposition \ref{dual_real} (ii),
\[\bigcup_{(d_0,\ldots,d_{m+L}) \in S^+} A_M(d_0,\ldots,d_{m+L}) \setminus \bigcup_{(d_0,\ldots,d_{m+L}) \in S^-} A_M(d_0,\ldots,d_{m+L})\]
is an (inhomogeneous) Bohr set $\mathcal{B}_{\alpha}(x,J)$ with $\lambda(J) \ll \gamma(\tfrac{2}{3})^L$.
An application of Lemma \ref{bohr_card} now proves \eqref{sandwich_error}.
\end{proof}

Proposition \ref{cylinder_cover} reduces Theorem \ref{thm_bohr_div} to establishing equidistribution on $\mathbb{Z}_d$ on Ostrowski cylinder sets $A_M(d_1,\ldots,d_m)$. This will be finally done in the subsequent technical key lemma.

\begin{lem}\label{markov_lem}
    Let $\alpha$ be badly approximable.
    For $x \in \N$, let $k = k(x)$ be such that 
    $q_{k-1} \leq x < q_{k+1}$ 
    and let $m = m(x) \in \N$ be such that $\lim_{x\to \infty}k(x) - m(x) = \infty$.
    Further, let $L$ be a fixed integer, let
    \[S_L = \left\{x' \in \mathbb{N}: x' = \sum_{i = k(x)-L+1}^{k(x)}b_iq_i, \quad (b_i)_{i} \text{ admissible }\right\},\] and let $(x_j)_{j \in \N}$
    be the sequence that enumerates the elements of $S_L$ in ascending order.
    Then for every $d \in \mathbb{N}$, every $a \in \mathbb{Z}_d$, and any admissible sequence $(d_i)_{i \in \N}$,
  \begin{equation}\label{bohr_div_special}
    \#\{N \in A_{x_j}(d_0,\ldots,d_m): N \equiv a \pmod d\} \sim \frac{\#A_{x_j}(d_0,\ldots,d_m)}{d},\quad  j \to \infty,
\end{equation}
with the convergence being uniform in $(d_i)_{i \in \N}$.
\end{lem}

\begin{proof}
We endow the finite set $A_{x_j}(d_0,\ldots,d_m)$ with the uniform probability measure $\mathbb{P} = \mathbb{P}_j$. By the Ostrowski expansion $N = \sum_{i \leq m(x_{j})}d_iq_i + \sum_{m(x_{j})+1 \leq i \leq k(x_{j})}b_iq_{i}$ we have the correspondence

\[N \longleftrightarrow (b_{m(x_{j})+1},\ldots,b_{k(x_{j})}),\]
so we can interpret the Ostrowski digits $b_i$ as random variables on $\{0,\ldots,a_{i+1}\}$. Note that they are neither uniformly distributed nor independent - a fact that arises from the Ostrowski rule $b_i = a_{i+1} \implies b_{i-1} = 0$.\\

We now define the set of \textit{free indices}\footnote{The term is motivated by the fact that the first $m$ digits are fixed and there is also a restriction for the last $L$  (by the Ostrowksi rule this might also restrict the choices at positions $m+1$ and $k(x_{j})-L-1$), but the probability distribution for $b_i$ is non-degenerated for $i \in F$.}, denoted by 
$F := F_{j,L} := \{m(x_{j})+2,\ldots,k(x_{j})-L-1\}$. 
Note that by the assumption of $\lim_{x \to \infty}k(x) - m(x) = \infty$, it follows that $\lim_{j \to \infty}\#F_{j,L} = \infty.$\\

We now claim to have the following three crucial facts for all $i \in F$:

\begin{itemize}
    \item[(i)] We have a uniform lower bound for the event of $b_i = 0$, even after fixing an arbitrary amount of other digits. Put formally, the random variable $\mathbb{P}[b_i=0|(b_1,\ldots,b_{i-1},b_{i+1},\ldots)]$ satisfies with probability $1$

    \[\mathbb{P}[b_i=0|(b_1,\ldots,b_{i-1},b_{i+1},\ldots)] \geq \frac{1}{C+1}.\]
    \item[(ii)] Let $Z_i := (b_i,b_{i+1})|b_{i+2} = 0, b_{i-1}= 0$, i.e., the random variable of $(b_i,b_{i+1})$ given the event $b_{i+2} = 0, b_{i-1}= 0$. Further, let $S_i$ be the set of all admissible tuples for $(b_i,b_{i+1})$, i.e.,   
    \[S_i:= \{(x,y) \in \mathbb{Z}: 0 \leq x \leq a_{i+1}, 0 \leq y \leq a_{i+2}-1\} \cup \{a_{i+2},0\}.\] Then for all $(x,y) \in S_i$,
    $\mathbb{P}[Z_i = (x,y)] = \frac{1}{|S_i|}$.
    \item[(iii)] $Z_i$ is independent of all digits $b_j$, $j \notin \{i-1,i,i+1,i+2\}$, i.e. formally:
    We have (almost surely)
    \[\mathbb{P}[Z_i = (x,y)|(b_1,\ldots,b_{i-2},b_{i+2},\ldots)] = \mathbb{P}[Z_i = (x,y)] = \frac{1}{|S_i|}, \quad (x,y) \in S_i.\]
\end{itemize}
Indeed, we observe that if $(d_0,\ldots,d_m,b_{m+1},\ldots,b,\ldots,b_{k(x_j)})$
 is admissible, then so is $(d_0,\ldots,d_m,b_{m+1},\ldots,0,\ldots,b_{k(x_j)})$.
Hence we find for all $b \in \{0,1,\ldots,a_{i+1}\}, i \in F$ an injection
\[\begin{split}\{N \in A_{x_j}(d_0,\ldots,d_m): b_i(N) = b\} &\to \{N \in A_{x_j}(d_0,\ldots,d_m): b_i(N) = 0\}\\N &\mapsto N - b_iq_i,
\end{split}\]
and this property remains valid by considering restrictions on both sets by fixing an arbitrary amount of digits $b_j, j \neq i$, proving (i).\\

By the same arguments, we find for all $(x,y),(x',y') \in S_i, i \in F$ a bijection between
\[\begin{split}&\{N \in A_{x_j}(d_0,\ldots,d_m): (b_i(N),b_{i+1}(N)) = (x,y), b_{i-1}(N)=  b_{i+1}(N) = 0\} \text{ and} \\
&\{N \in A_{x_j}(d_0,\ldots,d_m): (b_i(N),b_{i+1}(N)) = (x',y'), b_{i-1}(N)=  b_{i+1}(N) = 0\} 
\end{split}\]
via $N \mapsto N - (xq_{i+1} + y{q_i}) + (x'q_{i+1} + y'{q_i})$.
We remark that here the condition on $b_{i-1}(N)=  b_{i+1}(N) = 0$, as well as the special form of $x_j$ is used; without that, this statement would not be true in general.
For (iii), we again restrict both sets by fixing digits, still keeping the bijectivity property.\\

We now define $F_4 := \{i \in \N: \{4i,4i+1,4i+2,4i+3\} \in F\}$ and define the random variables
\[X_i := \mathds{1}_{[b_{4i} = b_{4i+3} = 0]}, \quad i \in F_4.\]
We observe that by (i), we have
\begin{equation}\label{lower_bound_0}\mathbb{P}[X_i = 1|E] > \frac{1}{(C+1)^2} \quad \text{ for any } E \in \sigma(\{b_j, j \notin \{4i,4i+3\}\}).\end{equation}
Next, we define
\[G_4 := \{4i \in \N: i \in F_4\} \cup \{4i+3 \in \N: i \in F_4\}, \quad H_4 = \{(h_i)_{i \in G_4}: (h_i)_{i \in G_4} \text{ is admissible }\}.\]
This allows us to partition
\[A_{x_j}(d_0,\ldots,d_m) = \dot{\bigcup_{\bm{h} \in H_4}}A_{\bm{h}},\]
where 
\[A_{\bm{h}} = \{N \in A_{x_j}(d_0,\ldots,d_m): b_i(N) = h_i \,\forall i \in G_4\}.\]

Clearly, $X_i$ is deterministic on $A_{\bm{h}}$ for all $i \in F_4$. We now focus on those $\bm{h}$
where $\sum_{i \in F_4}X_i \geq \log \log |F_4| =: R$ since this is a generic property:
Indeed, we have by using \eqref{lower_bound_0} that
\begin{equation}\label{slln}\mathbb{P}\left[\sum_{i\in F_4} X_i \leq R\right]
\leq \sum_{\substack{J \subseteq F_4\\ |J| = R}}\mathbb{P}\left[\sum_{i\in I\setminus J} X_i = 0\right]
\leq |F_4|^{R} \left(1 - \frac{1}{(C+1)^2}\right)^{|F_4|/2} \to 0, \quad j \to \infty.
\end{equation}
Thus we can assume from now on that $\bm{h}$ is of the form such that $\sum_{i \in F_4}X_i \geq R$.
This allows us to find for fixed $\bm{h}$ an index set $I = I_{\bm{h}} \subseteq F_4$ with $|I| = R$
and 
\[X_i = 1 \quad \forall i \in I.\]

By multiple applications of the pigeonhole principle, we obtain the existence of integers
$1 \leq c,c' \leq C(\alpha)$ such that there are 
$R'$ many $i \in I$ with $R' \gg_{d,C(\alpha)} R$ such that:
\begin{itemize}
    \item $q_{4i+1} \equiv a \pmod d,$
    \item $a_{4i+2} = c$.
    \item $q_{4i+2} \equiv b \pmod {d}$,
    \item $a_{4i+3} = c'$.
\end{itemize}
Since successive denominators of convergents are always coprime, we observe that $\gcd(a,b,d) = 1$.
Calling the index set where the above conditions hold $J_{\bm{h}} \subset I_{\bm{h}}$, 
we note that
$J_{\bm{h}} \gg_{d,\alpha} \log \log |F_4|$ and is thus unbounded. 

We will now prove equidistribution $\mod d$ on this set $A_{\bm{h}}$, with the convergence rate being uniform in $\bm{h}$. In order to do so, we define $\mathbb{P}_{\bm{h}}$ as the uniform measure on $A_{\bm{h}}$. Writing
\[N = \sum_{i \in J_{\bm{h}}}b_iq_i + \sum_{i \notin J_{\bm{h}}}b_iq_i,\]
we let $Y_0(N) = Y_0^{(\bm{h})}(N) := \sum_{i \notin J_{\bm{h}}}b_iq_i  \pmod d$ and let $\nu_{\bm{h}}$ be the distribution vector of $Y_0$ (with respect to $\mathbb{P}_{\bm{h}}$).
Since for all $i \in J_{\bm{h}}$ we have 
$q_{4i+1} \equiv a \pmod d, q_{4i+2} \equiv b \pmod d$
as well as $b_{4i} = b_{4i+3} = 0$, we get
\[N \equiv Y_0(N) + \sum_{i \in J_{\bm{h}}}ab_{4i+1}(N) + bb_{4i+2}(N) \pmod d.\]

Now setting \[Y_i(N) = Y_i^{(\bm{h})}(N) := ab_{4i+1}(N) + bb_{4i+2}(N) \pmod d, \quad i \in J_{\bm{h}},\]
we see that $Y_i$ is a random variable on $A_{\bm{h}}$ that is determined by $(b_{4i+1},b_{4i+2})$. Writing $S := S_{4i_0+1}$ for an arbitrary $i_0 \in J_{\bm{h}}$, we get by (ii), (iii), and our construction of $J_{\bm{h}}$ that 
\begin{equation}\label{ud_ph}\mathbb{P}_{\bm{h}}[(b_{4i+1},b_{4i+2})=(x,y)] = \frac{1}{|S|}\end{equation}
for all $(x,y) \in S$ and all $i \in J_{\bm{h}}$. Note that $(b_{4i+1},b_{4i+2})$ is (with respect to $\mathbb{P}_{\bm{h}}$) independent of
$\{(b_{4j+1},b_{4j+2})\}_{j \in J_{\bm{h}}\setminus \{i\}}$, a fact arising from (ii). Thus $(Y_i)_{i \in {J_{\bm{h}}}}$ is a sequence of i.i.d. random variables satisfying the following:

\begin{itemize}
    \item[(I)] $\mathbb{P}_{\bm{h}}[Y_i = j] = p_j^{(\bm{h})} =  \frac{u_j}{v_j}\in \mathbb{Q}$ with $v_j \ll_{C(\alpha),d} 1, \quad 0 \leq j \leq d-1$.
    \item[(II)] $\mathbb{P}_{\bm{h}}[Y_i = a] > 0, \mathbb{P}_{\bm{h}}[Y_i = b] > 0, \mathbb{P}_{\bm{h}}[Y_i = 0] > 0$.
    \item[(III)] $(Y_i)_{i \in J_{\bm{h}}}$ is independent of $Y_0$.
\end{itemize}
Indeed, (I) follows from $|S| \ll_{\alpha} 1$ and \eqref{ud_ph}, (II) from the fact that $\{(0,1),(0,0),(1,0)\} \subseteq A_i$, and (III) follows from (iii).

By relabeling, we get random variables $(X_i)_{0 \leq i \leq |J_{\bm{h}}|}$
with $X_0 \stackrel{d}{=} \theta_{\bm{h}}$, and $X_i \stackrel{d}{=} (p_1^{\bm{(h)}},\ldots,p_d^{\bm{(h)}})$. 

Thus by applying Corollary \ref{rand_walk}, we obtain that for any $\varepsilon > 0$
there exists $K_0((p_1^{\bm{(h)}},\ldots,p_d^{\bm{(h)}}),\varepsilon)$ such that
\begin{equation}\label{applic_conv}\left\lvert\mathbb{P}_{\bm{h}}\left[Y_0 + \sum_{i \in J_{\bm{h}}} Y_i \equiv a \mod d\right] -\frac{1}{d}\right\rvert < \varepsilon,\end{equation}
provided that $|J_{\bm{h}}| \geq K_0$. We stress that $K_0$ only depends on $(p_1^{\bm{(h)}},\ldots,p_d^{\bm{(h)}})$, and not on the initial distribution $\theta_{\bm{h}}$, and by (I),
there are only a bounded number number of possibilities for the vector $(p_1^{\bm{(h)}},\ldots,p_d^{\bm{(h)}})$.
Hence we obtain uniformity in \eqref{applic_conv} among all $\bm{h}$ that satisfy $\sum_{i \in F_4}X_i \geq K_0$, by choosing $K_0$ maximal among the finite number of distributions.
Since by \eqref{slln}, the contribution from those $\bm{h}$ with $\sum_{i \in F_4}X_i < R$ is negligible, summing over all $\bm{h}$
proves the statement since $j \to \infty$ implies $R \to \infty$.
\end{proof}

\begin{proof}[Proof of Theorem \ref{thm_bohr_div}]
    Let $m = m(x)$ such that $|\delta_{m-1}| < \rho_x = \lambda(I_x) \leq |\delta_m|$
    and $n(x)$ such that $q_{n} \leq x < q_{n+1}$. Note that by assumption of $x \rho_x \to \infty$ and $\alpha$ badly approximable, this means that $n(x) - m(x) \to \infty.$
    Let $L \in \N$ be a fixed integer such that $2L < n(x) - m(x)$.

 Taking $x^{-}_L \leq x \leq x^+_L$ defined as in Proposition \ref{prop_N_pm}, we obtain
\begin{equation}\label{N_pm}x^-_L \leq x \leq x_L^+, \quad \frac{x^+_L}{x^-_L} = 1 + O((\tfrac{2}{3})^L).\end{equation}
By Proposition \ref{cylinder_cover}, we thus get

\[\bigcup_{(d_0,\ldots,d_{m+L}) \in S^-} A_{x^{-}_L}(d_0,\ldots,d_{m+L}) \subseteq B_{\alpha}(x,I_x) \subseteq
\bigcup_{(d_0,\ldots,d_{m+L}) \in S^+} A_{x_L^+}(d_0,\ldots,d_{m+L})\]
with $|S^{\pm}| \ll L^{C(\alpha)+1}$, and (with $\Delta$ denoting the symmetric difference),
\begin{equation}
    \label{sandwich_err_thm}
\left|\left(\bigcup_{(d_0,\ldots,d_{m+L}) \in S^{\pm}} A_{x^{\pm}_L}(d_0,\ldots,d_{m+L})\right) \Delta B_{\alpha}(x,I_x)
\right| \ll_\alpha {x\rho_x\left(\tfrac{3}{2}\right)^{-L}}.\end{equation}
Note that $A_{x^+_L}(d_0,\ldots,d_{m+L})$ satisfies all conditions of Lemma \ref{markov_lem} (with $x^+_L = x_j$, $m = m+L, k = n, \gamma = \rho_x$), which implies for all $(d_0,\ldots,d_{m+L}) \in S^+$,

\[\#\{N \in A_{x_L^+}(d_0,\ldots,d_{m+L}): d \equiv a \pmod d\} = \frac{\#A_{x^{+}_L}(d_0,\ldots,d_{m+L})}{d}(1 + o(1)),\quad  n(x) \to \infty.\]
By $x \to \infty$, this implies together with \eqref{N_pm}, \eqref{sandwich_err_thm} and Lemma \ref{bohr_card}

\[\begin{split}&\limsup_{x \to \infty}\frac{\#\{n \in \mathcal{B}(x,I_x): n \equiv a \pmod d\}}{\#\mathcal{B}(x,I_x)}
\\\leq &\;\limsup_{x \to \infty}\frac{\sum_{(d_0,\ldots,d_{m+L}) \in S^+}\#\{n \in A_{x_L^+}(d_0,\ldots,d_{m+L}): n \equiv a \pmod d\}}{x\rho_x(1 + o(1))}
\\=&\;\frac{1}{d} + O((\tfrac{3}{2})^{-L}).
\end{split}\]
In the same way, using $S^-$ and $x_L^-$ in place of $S^+$ and $N_L^+$, we get
\[\liminf_{N \to \infty}\frac{\#\{n \in \mathcal{B}(x,I_x): n \equiv a \pmod d\}}{\#\mathcal{B}(x,I_x)}
\geq \frac{1}{d} - O((\tfrac{3}{2})^{-L}).
\]
With $L \to \infty$, this concludes the proof of the theorem.
\end{proof}

\section{Proof of Theorem \ref{main_thm}}
\subsection{Averaging multiplicative functions over Bohr sets}\label{sec_ant_part}
In this part, we establish a version of \eqref{decor} by applying Theorem \ref{thm_bohr_div}, and various maneuvers from analytic number theory. This will be the  most crucial step.
We start with a tail estimate that will be useful in the proof of the main statement in this section (Lemma \ref{key_lem}).

\begin{lemma}[Anatomy of integers on Bohr sets]\label{anatomy_lem}
Let $\alpha$ be badly approximable, let $\mathcal{B} =\mathcal{B}_{\alpha}(x,\rho)$ be a Bohr set with cardinality $B \gg (\log \log x)^{10}$
and let $r \in \mathbb{R}_{>0}$ with $r < B^{1/3}$.
Then for all $h' \in [0,x]$ we have
\begin{equation}\label{anatomic_cond}\frac{1}{B}\cdot\#\left\{h \in B: \prod_{\substack{p \mid h-h'\\p > r}}\left(1 + \frac{1}{p}\right)\geq 1 + \frac{1}{\sqrt{r}}\right\}  \ll \frac{1}{r^{1/3}}.\end{equation}
\end{lemma}

\begin{proof}
    By straightforward Taylor estimates, it suffices to show 
    \[\frac{1}{B}\#\left\{h \in \mathcal{B}: \sum_{\substack{p \mid h-h'\\p > r}}\frac{1}{p}\geq \frac{1}{2\sqrt{r}}\right\}  \ll \frac{1}{r^{1/3}}.\]
    Applying Markov's inequality (after removing if necessary the $h'$ from $\mathcal{B}$ which is clearly negligible) yields

    \[\frac{1}{B}\#\left\{h \in \mathcal{B}: \sum_{\substack{p \mid h-h'\\p > r}}\frac{1}{p}\geq \frac{1}{2\sqrt{r}}\right\}  \ll \frac{\sqrt{r}}{B}\sum_{h \in \mathcal{B}}\sum_{\substack{h - h' \equiv 0 \pmod p\\p > r}}\frac{1}{p}
   = \frac{\sqrt{r}}{B}\sum_{\substack{p > r}}\frac{1}{p}\sum_{\substack{h \in \mathcal{B}\\h \equiv h' \pmod p}}1.
    \]
Applying Proposition \ref{HK_machine}, we get 
\[\frac{1}{p}\sum_{\substack{h \in \mathcal{B}\\h \equiv h' \pmod p}}1
\ll_{\alpha} \frac{B}{p} + O(\sqrt{B}),
\]
which we apply for $p < x$, whereas for $p > x$, since $\mathcal{B} \subseteq [1,x]$, we trivially have for all $h \neq h'$
\[\frac{1}{p}\sum_{\substack{h \in \mathcal{B}\\h \equiv h' \pmod p}}1 = 0.\]
Thus we obtain

    \[\frac{1}{B}\#\left\{h \in \mathcal{B}: \sum_{\substack{p \mid h-h'\\p > r}}\frac{1}{p}\geq \frac{1}{2\sqrt{r}}\right\}  \ll \frac{1}{B} + \sqrt{r}\sum_{p > r}\frac{1}{p^2} + 
\frac{\sqrt{r}}{\sqrt{B}}\sum_{r < p < x}\frac{1}{p}
\ll \frac{1}{\sqrt{r}} + \frac{\sqrt{r}\log \log x}{\sqrt{B}} \ll \frac{1}{r^{1/3}},
    \]
    where we used $B \gg (\log \log x)^{10}$ and $r < B^{1/3}$ in the last step.
\end{proof}

\begin{lem}[Key Lemma]\label{key_lem}
    Let $k \geq 2$ be a fixed integer, $\bm{R} = [0,\pm u_1]\times \ldots \times [0,\pm u_{k-1}], u_i \in \mathbb{R}_{> 0}$
    and let \[\bm{B} = \bm{B}(\bm{R},N,\alpha) = \pm \mathcal{B}_{\alpha}(x_1,\pm u_1/N)\times \ldots \times \pm \mathcal{B}_{\alpha}(x_{k-1},\pm u_{k-1}/N)\] where $x_i = x_i(N), x = x(N)$ satisfies
    \[N^{10} \gg x \geq \max_{1 \leq i \leq k-1}x_i \geq \min_{1 \leq i \leq k-1} x_i \gg_A N (\log \log  N)^A\] for all $A > 0$.
    Then 
    \[\sum_{\bm{h} \in \bm{B}_{\neq}} \prod_{p \leq f(x)}\left(1 - \frac{g_{\bm{h}}(p)}{p}\right)
    \sim \#\bm{B} \prod_{p \leq f(x)} \left(1 - \frac{1}{p}\right)^k, \quad N \to \infty.
    \]
\end{lem}

\begin{proof}
We may assume that $\bm{B} = \bm{B}(\bm{R},N,\alpha) =  \mathcal{B}_{\alpha}(x_1,\pm u_1/N)\times \ldots \times \mathcal{B}_{\alpha}(x_{k-1},\pm u_{k-1}/N)$; the cases where some of the coordinates have negative entries only need obvious modifications.
    In the first step, we do a pre-sieving for primes $p \leq k$. For given $\bm{\ell} \in \mathbb{Z}_{k!}^{k-1}$, we define
    \[\bm{\mathcal{B}^{(\ell)}} := \{\bm{h} \in \bm{B}: \bm{h} \equiv \bm{\ell} \pmod{k!}\}.\]
Since $k$ is fixed throughout the proof, we remark for later reference that an application of Theorem \ref{thm_bohr_div} implies 
\begin{equation}\label{presieve}\#\bm{\mathcal{B}^{(\ell)}} \sim \frac{\#\bm{B}}{k!^{k-1}}.\end{equation}

Now we rewrite
    \begin{equation}\label{rewrite_fg}\begin{split}&\sum_{\bm{h} \in \bm{B}_{\neq}} \prod_{p \leq f(x)}\left(1 - \frac{g_{\bm{h}}(p)}{p}\right)
    = 
    \\&\sum_{\bm{\ell} \in \mathbb{Z}_{k!}^{k-1}} \prod_{1 \leq p \leq k} \left(1 - \frac{g_{\bm{\ell}}(p)}{p}\right)
    \prod_{k < p \leq f(x)}\left(1 - \frac{1}{p}\right)^{k} \prod_{k < p \leq f(x)}\frac{\left(1 - \frac{k}{p}\right)}{\left(1 - \frac{1}{p}\right)^{k}
}\sum_{\bm{h} \in \bm{\mathcal{B}^{(\ell)}}_{\neq}}\prod_{k < p \leq f(x)} f_{\bm{h}}(p),
\end{split}
    \end{equation}
where $f_{\bm{h}}$ is a multiplicative function that satisfies
\begin{equation}\label{property_f}1 \leq f_{\bm{h}}(p) \leq 1 + \frac{k - g_{\bm{h}}(p)}{p} + O(1/p^2), \quad k < p \leq f(x).\end{equation}

Next, we remove the contribution of large primes. We let $r_1 = (\log \log x)^C$ for suitably large $C = C(k)$ determined later. 
We claim to have the upper-bound
\begin{equation}\label{crude_f_bound}
f_{\bm{h}}(p) \leq \left(1 + O\left(\frac{1}{p^2}\right)\right)\prod_{i = 1}^{k-1} \prod_{p \mid h_i}\left(1 + \frac{k}{p}\right)
\prod_{i = 1}^{k-1}\prod_{j = i+1}^{k-1}\prod_{p \mid h_i - h_j}\left(1 + \frac{k}{p}\right) =: \tilde{f}_{\bm{h}}(p).
\end{equation}
Indeed, if $p \nmid h_i$ for all $i$ and $p \nmid h_i - h_j \forall i \neq j$, then 
$g_{\bm{h}}(p) = \#\{0^{(p)},h_{1}^{(p)},\ldots,h_{k-1}^{(p)}\} = k$ and thus $f_{\bm{h}}(p) \leq 1 + O(\frac{1}{p^2})$ by \eqref{property_f}. In all other cases, we simply use $g_{\bm{h}}(p) \geq 0$ in \eqref{property_f}, which proves \eqref{crude_f_bound}.
Assuming for a fixed $\bm{h} \in \bm{B}_{\neq}$ and $x$ sufficiently large that
\[\prod_{r_1 \leq p \leq f(x)}f_{\bm{h}}(p) \geq 1 + \frac{1}{r_1^{1/3}},\]
there must exist $1 \leq i \leq k-1$ or $1\leq i < j\leq k-1$
such that

\[\prod_{p \mid h_i}\left(1 + \frac{k}{p}\right) \geq 1 + \frac{1}{r_1^{3/5}} \quad \text{ or } \quad \prod_{p \mid h_i-h_j}\left(1 + \frac{k}{p}\right) \geq 1 + \frac{1}{r_1^{3/5}},\]
which implies 
\begin{equation}\label{anatomic_prop}\prod_{p \mid h_i}\left(1 + \frac{1}{p}\right) \geq 1 + \frac{1}{r_1^{1/2}} \quad \text{ or } \quad \prod_{p \mid h_i-h_j}\left(1 + \frac{1}{p}\right) \geq \frac{1}{r_1^{1/2}}.\end{equation}

 An application of Lemma \ref{anatomy_lem} together with \eqref{presieve} now implies that any fixed case in \eqref{anatomic_prop} happens in $\bm{B}^{(\bm{\ell})}$ for at most $\ll_k \frac{\bm{B}^{(\bm{\ell})}}{r_1^{1/3}}$ many elements. By employing the union bound, we deduce
\begin{equation}\label{anatomic_k_tuple}\frac{1}{\# \bm{B}^{(\bm{\ell})} }\#\left\{
\bm{h} \in \bm{B}^{(\bm{\ell})}: \prod_{p \geq r_1}f_{\bm{h}}(p) \geq 1 + \frac{1}{r_1^{1/3}}
\right\} \ll_k \frac{1}{r_1^{1/3}}.\end{equation}

Let $\bm{F}_1 \subset \bm{B}^{(\bm{\ell})}$ denote the set of those $\bm{h} \in \bm{B}^{(\bm{\ell})}_{\neq}$ where $\prod_{p \geq r_1}f_{\bm{h}}(p) \geq 1 + \frac{1}{r_1^{1/3}}$ and set $\bm{E}_1 := \bm{B}^{(\bm{\ell})}_{\neq}\setminus \bm{F}_1$.
We start to estimate the contribution of $\bm{F}_1$: Note that 
by \eqref{crude_f_bound}, we obtain that for any $\bm{h}\in \bm{B}_{\neq}$
\begin{equation}\label{pointwise_f_bound}\begin{split}\prod_{k < p < f(x)}f_{\bm{h}}(p) &\ll \prod_{i = 1}^{k-1} \prod_{p \mid h_i}\left(1 + \frac{1}{p}\right)^k \prod_{i = 1}^{k-1}\prod_{j = i+1}^{k-1}\prod_{p \mid h_i - h_j}\left(1 + \frac{1}{p}\right)^k \\&\ll_k \prod_{i = 1}^{k-1}\left(\frac{h_i}{\varphi(h_i)}\right)^k \prod_{i = 1}^{k-1}\prod_{j = i+1}^{k-1}\left(\frac{|h_i-h_j|}{\varphi(|h_i-h_j|)}\right)^k 
\ll_k (\log \log x)^{k^3},\end{split}\end{equation}
by the classical upper bound $\frac{m}{\varphi(m)} \ll \log \log m$. 
In combination with \eqref{anatomic_k_tuple}, this shows

\[\frac{1}{\# \bm{B}^{(\bm{\ell})}}\sum_{\bm{h} \in \bm{F}_1} \prod_{k < p < f(x)}f_{\bm{h}}(p) \ll_k \frac{(\log \log x)^{k^3}}{(\log \log x)^{C(k)/3}} = o(1),\]
when choosing $C >3k^3$.
For $\bm{E}_1$, we have by definition that
\[1 \leq \prod_{p \geq r_1}f_{\bm{h}}(p) \leq 1 + \frac{1}{r_1^{1/3}},\]
and an analogous calculation to the above shows 

\[\frac{1}{\# \bm{B}^{(\bm{\ell})}}\frac{1}{r_1^{1/3}}\sum_{\bm{h} \in \bm{E}_1} \prod_{k < p < r_1}f_{\bm{h}}(p) \ll_k \frac{(\log \log x)^{k^3}}{(\log \log x)^{C(k)/3}} = o(1),\]
since $C >3k^3$.
Thus we have shown that
\begin{equation*}\sum_{\bm{h} \in \bm{B}^{(\bm{\ell})}_{\neq}}\prod_{k < p < f(x)}f_{\bm{h}}(p) \sim \sum_{\bm{h} \in \bm{E}_1}\prod_{k < p < r_1}f_{\bm{h}}(p),\end{equation*}
where $r_1 = (\log \log x)^{3k^3 + 1}$.\\

We now define $r_2 = (\log \log \log x)^{C'}$ for some $C' = C'(k)$ chosen later. Applying Lemma \ref{anatomy_lem}, now to $r_2$, and arguing as in \eqref{anatomic_prop} shows that the set
\[\bm{F}_2 := \left\{\bm{h} \in \bm{E}_1: \prod_{p > r_2}f_{\bm{h}}(p) \geq 1 + \frac{1}{r_2^{1/3}}\right\}\]
satisfies 
\[\frac{\# \bm{F}_2}{\# \bm{B}^{(\bm{\ell})} } \ll_k \frac{1}{r_2^{1/3}}.\]
Note that on $\bm{F}_2 \subseteq \bm{E}_1$, we have the sharper estimate 

\[\prod_{k < p < f(x)}f_{\bm{h}}(p) \ll \prod_{k < p < r_1}f_{\bm{h}}(p)
\ll_k \left(\prod_{p \leq r_1}\left(1 + \frac{1}{p}\right)\right)^{k^3}
\ll (\log r_1)^{k^3} \ll_k (\log \log \log x)^{k^3} = o(r_2^{1/3}),
\]
when choosing $C' > 3k^3$. For $\bm{E}_1 \setminus \bm{F}_2$, we again argue as before, showing that the error term $1/r_2^{1/3}$ is also negligible.
Thus we have removed the contribution of all primes $p > r_2$, i.e. we obtain 
\begin{equation}\label{contr_large}\sum_{\bm{h} \in \bm{B}^{(\bm{\ell})}_{\neq}}\prod_{k < p < f(x)}f_{\bm{h}}(p) \sim \sum_{\bm{h} \in \bm{B}^{(\bm{\ell})}_{\neq}}\prod_{k < p < r_2}f_{\bm{h}}(p),\end{equation}
where $r_2 = (\log \log \log x)^{3k^3 + 1}$.

We now define $s_{\bm{h}}(p) := p(f_{\bm{h}}(p) - 1)$ for $p < k < r_2$ and extend this function multiplicatively via $s_{\bm{h}}(d) := \mu^2(d)\prod_{p \mid d}s_{\bm{h}}(p)$.
This allows us to write

\[\prod_{k < p \leq r_2} f_{\bm{h}}(p) = \sum_{d \mid P(r_2)}\frac{\mu^2(d)s_{\bm{h}}(d)}{d}\]
where $P(r_2) = \prod_{k < p < r_2}p$.\\

Note that $s_{\bm{h}}(d)$ only depends on the residue classes of $\bm{h} \pmod d$. Thus we have
\[\sum_{\bm{h} \in \bm{B}^{(\bm{\ell})}_{\neq}}\prod_{k < p < r_2}f_{\bm{h}}(p) = \sum_{d \mid P(r_2)} \frac{\mu^2(d)}{d}\sum_{\bm{j} \in \mathbb{Z}_d^{k-1}} s_{\bm{j}}(d)\sum_{\substack{\bm{h} \in \bm{B}^{(\bm{\ell})}_{\neq}\\\bm{h} \equiv \bm{j} \pmod d}}1.\]

Next, we fix a large integer $M$ that will be sent to $\infty$ after $N \to \infty$.
For $d < M$, we use equidistribution mod $d$ in the form of Theorem \ref{thm_bohr_div}:
Note that we have $\# \bm{B}_{\neq} \sim \#\bm{B}$ and therefore, by (note that $\gcd(P(r_2),k!) = 1$) Theorem \ref{thm_bohr_div} applied to every coordinate, we have for $d \mid P(r_2)$

\[\sum_{\substack{\bm{h} \in \bm{B}^{(\bm{\ell})}_{\neq}\\\bm{h} \equiv \bm{j} \pmod d}}1
= \sum_{\substack{\bm{h} \in \bm{B}^{(\bm{\ell})}\\\bm{h} \equiv \bm{j} \pmod d}}1 + o(\#\bm{B})
\sim \frac{\#\bm{B}}{d^{k-1}(k!)^{k-1}}.\]
Consequently,

\begin{equation}\label{small_d}\sum_{\substack{d \mid P(r_2)\\d \leq M}} \frac{\mu^2(d)}{d} \sum_{\bm{j} \in \mathbb{Z}_d^{k-1}}s_{\bm{j}}(d)\sum_{\substack{\bm{h} \in \bm{B}^{(\bm{\ell})}_{\neq}\\\bm{h} \equiv \bm{j} \pmod d}}1
\sim \frac{\# \bm{B}}{(k!)^{k-1}} \sum_{\substack{d \mid P(r_2)\\d \leq M}}\frac{\mu^2(d)}{d^k} \sum_{\bm{j} \in \mathbb{Z}_d^{k-1}}s_{\bm{j}}(d).
\end{equation}

Next, we will show that 

\begin{equation}\label{convergent_tail}\sum_{\substack{d \mid P(r_2)}}\frac{\mu^2(d)}{d^k} \sum_{\bm{j} \in \mathbb{Z}_d^{k-1}}s_{\bm{j}}(d) < \infty,\end{equation}
and consequently, 
\begin{equation}\label{vanishing_tail}\sum_{\substack{d \mid P(r_2)\\d > M}}\frac{\mu^2(d)}{d^k} \sum_{\bm{j} \in \mathbb{Z}_d^{k-1}}s_{\bm{j}}(d) = o(1), \quad M \to \infty.\end{equation}
To show \eqref{convergent_tail}, we prove that

\[\sum_{\bm{j} \in \mathbb{Z}_d^{k-1}}s_{\bm{j}}(d) \ll_{k,\varepsilon} d^{k-2 + \varepsilon}.\]
This is obviously trivial for non-squarefree $d$, so we can focus on squarefree $d$.
Note that from the multiplicative definition of $s_{\bm{j}}(d)$ and the Chinese Remainder Theorem, 
\begin{equation}\label{CRT}\sum_{\bm{j} \in \mathbb{Z}_d^{k-1}}s_{\bm{j}}(d) 
= \prod_{p \mid d} \sum_{\bm{j} \in \mathbb{Z}_p^{k-1}}s_{\bm{j}}(p).
\end{equation}
Using again \eqref{crude_f_bound}, we see that for $\bm{j} \in \mathbb{Z}_p^{k-1}$
with $p \nmid j_i$ and $p \nmid j_i - j_{\ell}$ for $1 \leq i < \ell \leq k-1$, we have
$f_{\bm{j}}(p) \leq 1 + O(p^{-2})$, so $s_{\bm{j}}(p) = p(f_{\bm{j}}(p)-1) \ll \frac{1}{p}$. Thus the contribution of these $\bm{j} \in \mathbb{Z}_p^{k-1}$ is $\ll \frac{p^{k-1}}{p}$. Note that there are at most $k^2 p^{k-2}$ many $\bm{j} \in \mathbb{Z}_p^{k-1}$ such that $p \mid j_i$ or $p \mid j_i - j_{\ell}$ for some $1 \leq i < \ell \leq k-1$, in which case \eqref{crude_f_bound} provides an upper bound of the form $s_{\bm{j}}(p)\ll \left(1 + \frac{1}{p}\right)^{k^3}$. Thus

\[\sum_{\bm{j} \in \mathbb{Z}_p^{k-1}}s_{\bm{j}}(p) \ll_k p^{k-2},\]
so \eqref{CRT} shows

\begin{equation}\label{same_coord_upper}\sum_{\bm{j} \in \mathbb{Z}_d^{k-1}}s_{\bm{j}}(d) \ll_k d^{k-2}O_k(1)^{\omega(d)}
\ll_{k,\varepsilon} d^{k-2 + \varepsilon}.
\end{equation}
Choosing $\varepsilon = 1/2$ finally shows \eqref{convergent_tail}, implying \eqref{vanishing_tail}.\\

Next, we will treat 
\[\sum_{\substack{d \mid P(r_2)\\d > M}} \frac{\mu^2(d)}{d} \sum_{\bm{j} \in \mathbb{Z}_d}s_{\bm{j}}(d)\sum_{\substack{\bm{h} \in \bm{B}^{\bm{(\ell)}}_{\neq}\\\bm{h} \equiv \bm{j} \pmod d}}1
\leq \sum_{\substack{d \mid P(r_2)\\d > M}} \frac{\mu^2(d)}{d} \sum_{\bm{j} \in \mathbb{Z}_d}s_{\bm{j}}(d)\sum_{\substack{\bm{h} \in \bm{B}\\\bm{h} \equiv \bm{j} \pmod d}}1,
\]
which we split into different ranges: Let $T = \sqrt{\min_{1 \leq i \leq k-1}(x_i/N)}$,
which by assumption satisfies $T \gg_A (\log \log N)^A$.
We apply Lemma \ref{bohr_card} to see that $\#\mathcal{B}(x_i,s_i) \gg_{s_i} T^2$.
Thus applying Proposition \ref{HK_machine} on each coordinate of $\bm{B}$, we obtain for any $M < d < T$ and any $\bm{j} \in \mathbb{Z}_d^{k-1}$,
\[\sum_{\substack{\bm{h} \in \bm{B}\\\bm{h} \equiv \bm{j} \pmod d}}1 \ll_{k} \frac{\#\bm{B}}{d^{k-1}}.\]
This shows 
\begin{equation}\label{interm_range}\sum_{\substack{d \mid P(r_2)\\M < d < T}} \frac{\mu^2(d)}{d} \sum_{\bm{j} \in \mathbb{Z}_d^{k-1}}\sum_{\substack{\bm{h} \in \bm{B}\\\bm{h} \equiv \bm{j} \pmod d}}s_{\bm{j}}(d) \ll_k \#\bm{B} \sum_{\substack{d \mid P(r_2)\\M < d < T}} \frac{\mu^2(d)}{d^k}\sum_{\bm{j} \in \mathbb{Z}_d^{k-1}}s_{\bm{j}}(d) = o(\#\bm{B}), \quad M \to \infty\end{equation}
by using \eqref{vanishing_tail} once more.

For $d > T$ we use \eqref{crude_f_bound} again to deduce  $s_{\bm{j}}(d) \ll \prod_{p \mid d}\left(1 + \frac{1}{p}\right)^{k^3} \ll \tau^{c_k}(d)$ for some $c_k$ only depending on $k$.
Furthermore, we observe the following immediate facts for squarefree numbers $d \mid P(r_2)$:

\begin{itemize}
    \item $d \leq r_2!$.
    \item Since $d$ is $r_2$-smooth, for $d > T$, we have
$T < r_2^{\omega(d)}$, hence
    \[\tau(d) = 2^{\omega(d)} \geq \exp\left(\frac{\log T}{\log r_2}\log 2\right) =: U.\]
\end{itemize}
Thus we get by using Rankin's trick

\[\begin{split}
    \sum_{\substack{d \mid P(r_2)\\d > T}}\frac{\mu^2(d)}{d}\sup_{\bm{j} \in \mathbb{Z}_d^{k-1}}s_{\bm{j}}(d)
    &\ll \sum_{\substack{d \mid P(r_2)\\ d > T}}\frac{\mu^2(d)\tau^{c_k}(d)}{d}
    \ll \frac{1}{U} \sum_{\substack{d \leq r_2!}}\frac{\mu^2(d)\tau^{c_k+1}(d)}{d}
    \ll \frac{(r_2 \log r_2)^{2^{c_k+1}}}{U},
\end{split}\]
where we used Stirling's formula and the well-known estimate
$\sum_{n \leq x}\frac{\tau(n)^{m}}{n} \ll_m (\log x)^{2^m}$.
Since by assumption $T \gg_A (\log \log x)^A$ for all $A > 0$ and $r_2 = (\log \log \log x)^{C'}$, we get that $U \gg_A (\log \log x)^A $ for all $A > 0$ and thus
$(r_2 \log r_2)^{2^{c_k}+1} = o(U)$.
Hence we obtain 
\begin{equation}\label{rankin_conc}\sum_{\substack{d \mid P(r_2)\\d>T}} \frac{\mu^2(d)}{d} \sum_{\bm{j} \in \mathbb{Z}_d^{k-1}}\sum_{\substack{\bm{h} \in \bm{B}\\\bm{h} \equiv \bm{j} \pmod d}}s_{\bm{j}}(d) = o(\#\bm{B}), \quad N \to \infty.\end{equation}

A combination of \eqref{contr_large}, \eqref{small_d}, \eqref{vanishing_tail}, \eqref{interm_range} and \eqref{rankin_conc} thus implies by $N \to \infty$ and then $M \to \infty$ that
\[\sum_{\bm{h} \in \bm{B}^{(\bm{\ell})}_{\neq}}\prod_{k < p < r_2}f_{\bm{h}}(p) \sim \frac{\#\bm{B}}{(k!)^{k-1}}\sum_{\substack{d \in \mathbb{N}\\\gcd(d,k!)= 1}}\frac{\mu^2(d)}{d^k} \sum_{\bm{j} \in \mathbb{Z}_d^{k-1}}s_{\bm{j}}(d).\]
Consequently, recalling \eqref{rewrite_fg}, we have

\begin{equation}\label{key_lem_final_compl}\begin{split}&\frac{1}{\#\bm{B}\prod_{p \leq f(x)} \left(1 - \frac{1}{p}\right)^k}
\sum_{\bm{h} \in \bm{B}_{\neq}} \prod_{p \leq f(x)}\left(1 - \frac{g_{\bm{h}}(p)}{p}\right)
    \\\sim &
    \sum_{\bm{\ell} \in \mathbb{Z}_{k!}^{k-1}} \left(\prod_{1 \leq p \leq k} \left(1 - \frac{g_{\bm{\ell}}(p)}{p}\right)\right)
    \prod_{p > k}\frac{\left(1 - \frac{k}{p}\right)}{\left(1 - \frac{1}{p}\right)^{k}
}\frac{1}{(k!)^{k-1}}\sum_{\substack{d \in \mathbb{N}\\\gcd(d,k!)= 1}}\frac{\mu^2(d)}{d^k} \sum_{\bm{j} \in \mathbb{Z}_d^{k-1}}s_{\bm{j}}(d),\end{split}\end{equation}
since $\frac{1 - \frac{k}{p}}{\left(1 - \frac{1}{p}\right)^{k}} = 1 + O_k(1/p^2)$ implies that
\[\prod_{k < p < f(x)}\frac{\left(1 - \frac{k}{p}\right)}{\left(1 - \frac{1}{p}\right)^{k}} \to \prod_{p > k}\frac{\left(1 - \frac{k}{p}\right)}{\left(1 - \frac{1}{p}\right)^{k}} < \infty, \quad x \to \infty.\]

We prove that the right-hand side of \eqref{key_lem_final_compl} evaluates to
$1$
by the following argument: Let 
$\bm{V} = [1,\prod_{p \leq f(x)}p]^{k-1}_{\neq}$. Since for all primes $p' \leq f(x), p' \mid \prod_{p \leq f(x)}p$, we can follow all the above steps (most steps without any error term at all)
to deduce
\[\begin{split}&\frac{1}{\prod_{p \leq f(x)} \left(1 - \frac{1}{p}\right)^k\#\bm{V}}\sum_{\bm{h} \in \bm{V}_{\neq}} \prod_{p \leq f(x)}\left(1 - \frac{g_{\bm{h}}(p)}{p}\right) \\\sim &\sum_{\bm{\ell} \in \mathbb{Z}_{k!}^{k-1}} \left(\prod_{1 \leq p \leq k} \left(1 - \frac{g_{\bm{\ell}}(p)}{p}\right)\right)\prod_{p > k}\frac{\left(1 - \frac{k}{p}\right)}{\left(1 - \frac{1}{p}\right)^{k}
} \frac{1}{(k!)^{k-1}} 
 \sum_{\substack{d \mid P(r_2)}}\frac{\mu^2(d)}{d^k} \sum_{\bm{j} \in \mathbb{Z}_d^{k-1}}s_{\bm{j}}(d).\end{split}\]
On the other hand, we use the Chinese Remainder Theorem to write
\[\frac{1}{\#\bm{V}}\sum_{\bm{h} \in \bm{V}} \prod_{p \leq f(x)}\left(1 - \frac{g_{\bm{h}}(p)}{p}\right)
= \prod_{p \leq f(x)}\frac{1}{p^{k-1}}\sum_{\bm{h} \in \mathbb{Z}_p^{k-1}}
\left(1 - \frac{g_{\bm{h}}(p)}{p}\right),
\]
claiming that
\begin{equation}\label{cute_id}\frac{1}{p^{k-1}}\sum_{\bm{h} \in \mathbb{Z}_p^{k-1}}
\left(1 - \frac{g_{\bm{h}}(p)}{p}\right) = \left(1 - \frac{1}{p}\right)^{k-1}.\end{equation}
Indeed, this follows from a simple induction procedure, noting that for $k \geq 2$,
\[\begin{split}\sum_{\bm{h} \in \mathbb{Z}_p^{k}}
\left(1 - \frac{g_{\bm{h}}(p)}{p}\right) &= 
\sum_{\substack{(\bm{h},j) \in \mathbb{Z}_p^{k}\\ g_{\bm{h}}(p) = g_{(\bm{h},j)}(p)}}\left(1 - \frac{g_{\bm{h}}(p)}{p}\right) + \sum_{\substack{(\bm{h},j) \in \mathbb{Z}_p^{k}\\ g_{\bm{h}}(p) = g_{(\bm{h},j)}(p) -1}}
\left(1 - \frac{g_{\bm{h}}(p)}{p}- \frac{1}{p}\right)
\\&= p\left(\sum_{\substack{\bm{h} \in \mathbb{Z}_p^{k-1}}}\frac{g_{\bm{h}}(p)}{p}\left(1 - \frac{g_{\bm{h}}(p)}{p}\right) + \sum_{\substack{\bm{h} \in \mathbb{Z}_p^{k-1}}}
\left(1 - \frac{g_{\bm{h}}(p)}{p}\right)\left(1 - \frac{g_{\bm{h}}(p)}{p}- \frac{1}{p}\right)\right)
\\&= p\left(\left(1 - \frac{1}{p}\right)\sum_{\substack{\bm{h} \in \mathbb{Z}_p^{k-1}}}\left(1 - \frac{g_{\bm{h}}(p)}{p}\right) 
\right)
\\&= p^k\left(1 - \frac{1}{p}\right)^{k},
\end{split}\]
where we used the induction hypothesis for $k-1$ in the very last line.
This proves \eqref{cute_id}, thus the proof of the lemma now follows from 
\[\frac{1}{\#\bm{V}}\sum_{\bm{h} \in \bm{V}\setminus \bm{V}_{\neq}} \prod_{p \leq f(x)}\left(1 - \frac{g_{\bm{h}}(p)}{p}\right)
\ll_k \frac{1}{\prod_{p \leq f(x)}p} = o\left({\prod_{p \leq f(x)} \left(1 - \frac{1}{p}\right)^k}\right).
\]
\end{proof}

\subsection{Step function procedure}\label{sec_step_fct}
Given an axis-parallel rectangle $\bm{R} \subset \mathbb{R}^{k-1}$, we may assume without loss of generality that
$\bm{R} = [0,\pm u_1]\times \ldots \times [0,\pm u_{k-1}]$ (where we understand $[a,b]$ for $b < a$ as $[b,a]$),
since we can otherwise partition $\bm{R}$ into $\ll_k 1$ many such rectangles.
We have by the very definition of the $k$-th order correlations

\[\begin{split}
 R_k^{(f)}(\bm{R},x) &= \frac{1}{\Phi(x,f(x))}\#\left\{\bm{n} = (\bm{m},n_k) \in [1,x]^{k-1}_{\neq}: P^{-}(\bm{n}) > f(x): \{\left(n_k\cdot \bm{1}_{k-1} - \bm{m}\right)\alpha\} \in \frac{\bm{R}}{\Phi(x,f(x))} \right\}
\\&= \frac{1}{\Phi(x,f(x))}\sum_{\bm{h} \in \bm{\hat{B}}_{\neq}}\#\left\{\bm{n} = (\bm{m},n_k) \in [1,x]^{k-1}: P^{-}(\bm{n}) > f(x): \left(n_k\cdot \bm{1}_{k-1} - \bm{m}\right) = \bm{h} \right\}
\end{split}
\]
where 
$\bm{\hat{B}} = \bm{\hat{B}}(\bm{R},x,\alpha) = \mathcal{\hat{B}}_{\alpha}(x,\pm \frac{u_1}{\Phi(x,f(x))})\times \ldots \times \mathcal{\hat{B}}_{\alpha}(x,\pm \frac{u_{k-1}}{\Phi(x,f(x))})$
with \[\mathcal{\hat{B}}_{\alpha}(x,\rho)  := \begin{cases}
    \{-x \leq n \leq x: \{n\alpha\} \in (0,\rho]\} \setminus \{0\} &\text{ for } \rho > 0,\\
    \{-x \leq n \leq x: \{n\alpha\} \in [1 - \rho,1)\} \setminus \{0\} &\text{ for } \rho < 0.
\end{cases}
\]
An application of Lemma \ref{sieve_k_tuple} shows

\begin{equation}\label{sieve_appl}\Phi(x,f(x))R_k^{(f)}(\bm{R},x) \sim \sum_{\bm{h} \in \bm{\hat{B}}_{\neq}}(x - (h^{+}_{>0}-h^{-}_{<0}))_{>0}
\prod_{p \leq f(x)}\left(1 - \frac{g_{\bm{h}}(p)}{p}\right) + O(x^{o(1)}).
\end{equation}

Next, we fix a large integer $T$ and let $I_T := \left\{\bm{t} = \left(\frac{j_1}{T},\ldots,\frac{j_{k-1}}{T}\right), j_i \in [-T,T)\cap\mathbb{Z} \right\}$. We say that an integer vector $\bm{h} \in [-x,x]^{k-1}$ satisfies
$\bm{h} \in S_{\bm{t}} = S_{\bm{t}}(x)$ if $h_i \in \left[\frac{xj_i}{T},\frac{x(j_i+1)}{T}\right]$ for all $1 \leq i \leq k-1$. In this notation, we observe that
\begin{equation}\begin{split}
&\sum_{\bm{h} \in \bm{\hat{B}}_{\neq}}(x - (h^+_{>0}-h^-_{<0}))_{> 0}
\prod_{p \leq f(x)}\left(1 - \frac{g_{\bm{h}}(p)}{p}\right)
\\=& \sum_{\bm{t} \in I_T}\sum_{\bm{h} \in \bm{\hat{B}}_{\neq} \cap S_{\bm{t}}}(x - (h^+_{>0}-h^-_{<0}))_{> 0}
\prod_{p \leq f(x)}\left(1 - \frac{g_{\bm{h}}(p)}{p}\right) 
\\=& \sum_{\bm{t} \in I_T}\sum_{\bm{h} \in \bm{\hat{B}}_{\neq} \cap S_{\bm{t}}}((x - (xt^+_{>0}-xt^-_{<0}))_{> 0}+ O(x/T))\prod_{p \leq f(x)}\left(1 - \frac{g_{\bm{h}}(p)}{p}\right)
\\=& x\sum_{\bm{t} \in I_T}(1 - (t^+_{>0}-t^-_{<0}))_{> 0}\sum_{\bm{h} \in \bm{\hat{B}}_{\neq} \cap S_{\bm{t}}}\prod_{p \leq f(x)}\left(1 - \frac{g_{\bm{h}}(p)}{p}\right) + O\left(\frac{x}{T}\sum_{\bm{h} \in \bm{\hat{B}}_{\neq}}\prod_{p \leq f(x)}\left(1 - \frac{g_{\bm{h}}(p)}{p}\right)  \right). \label{error_step_fct}
\end{split}
\end{equation}

Note that for every fixed $\bm{t} \in {I}_{T}$ we can write by simple inclusion-exclusion arguments $\bm{\hat{B}}_{\neq} \cap S_{\bm{t}}$ as a linear combination of $\leq 2^{k-1}$ many $\bm{B}_{\neq}$'s of the shape from Lemma \ref{key_lem}, with $\min_{1 \leq i \leq k-1}x_i \geq \frac{x}{T}$.
Thus we can apply Lemma \ref{key_lem} (since $x/T \gg_A ({\log \log N})^A$) which proves
that
\[\begin{split}&\;x\sum_{\bm{t} \in I_{T}}(1 - (t^+_{>0}-t^-_{<0}))_{> 0}\sum_{\bm{h} \in \bm{\hat{B}}_{\neq} \cap S_{\bm{t}}}\prod_{p \leq f(x)}\left(1 - \frac{g_{\bm{h}}(p)}{p}\right)  \\\sim&\;
x\prod_{p \leq f(x)} \left(1 - \frac{1}{p}\right)^k\sum_{\bm{t} \in I_{T}}(1 - (t^+_{>0}-t^-_{<0}))_{> 0}\#(\bm{\hat{B}} \cap S_{\bm{t}}).
\end{split}\]

Applying Lemma \ref{bohr_card} to $\bm{\hat{B}}_{\neq} \cap S_{\bm{t}}$ via another inclusion-exclusion argument proves
$\bm{\hat{B}}_{\neq} \cap S_{\bm{t}} \sim \frac{\#\bm{\hat{B}}}{\#I_T}$, which implies that the above is asymptotically equal to
\[\begin{split}
x^{k}\frac{\vol(\bm{R})}{\Phi(x,f(x))^{k-1}}\prod_{p \leq f(x)} \left(1 - \frac{1}{p}\right)^k\frac{1}{\#I_T}\sum_{\bm{t} \in I_{T}}(1 - (t^+_{>0}-t^-_{<0}))_{> 0}.
\end{split}\]
Another application of Lemma \ref{key_lem} and Lemma \ref{bohr_card} proves that the error term in \eqref{error_step_fct} 
is \[\ll \frac{x}{T}\sum_{\bm{h} \in \bm{\hat{B}}_{\neq}}\prod_{p \leq f(x)}\left(1 - \frac{g_{\bm{h}}(p)}{p}\right) \ll \frac{x^k}{T}\frac{\vol(\bm{R})}{\Phi(x,f(x))^{k-1}}\prod_{p \leq f(x)} \left(1 - \frac{1}{p}\right)^k,\]
and that the error term from \eqref{sieve_appl} is bounded above by 
\[O(\#\bm{\hat{B}} x^{o(1)}) = O_{\bm{R}}\left(\frac{x^{k-1 + o(1)}}{\Phi(x,f(x))^{k-1}}\right).\]
Since $f(x) = x^{o(1)}$, we have
$\Phi(x,f(x)) \sim x \prod\limits_{p \leq f(x)}\left(1 - \frac{1}{p}\right)$, and thus we get
\[\lim_{x \to \infty}R_k^{(f)}(\bm{R},x) = \vol(\bm{R})\left(\frac{1}{|I_T|}\sum_{\bm{t} \in I_T}(1 - (t^+_{>0}-t^-_{<0}))_{> 0} + O(1/T)\right).\]

We note that with $T\to \infty$, the above Riemann sum converges to its integral, i.e. 
\[\lim_{T \to \infty}\frac{1}{|I_T|}\sum_{\bm{t} \in I_T}(1 - (t^+_{>0}-t^-_{<0}))_{> 0} = \int_{[-1,1]^{k-1}}(1 - (y^+_{>0}-y^-_{<0}))_{> 0} \,\mathrm{d}\bm{y} = 1,\]
with the last equality following from a longer, but elementary calculation (see \cite[Lemma 6.2]{HZ_higher_order} for a rigorous proof). Hence we obtain

\[\lim_{x \to \infty}R_k^{(f)}(\bm{R},x) = \vol(\bm{R}),\]
as claimed. This proves the statement of Theorem \ref{main_thm} for the triangular arrays $((a_n^{(f)})_{n \leq x})_{x \in \mathbb{N}}$.

\subsection{From triangular arrays to sequences}\label{sec_arr_to_seq}

Let again w.l.o.g. $\bm{R} = [0,\pm u_1]\times \ldots \times [0,\pm u_{k-1}]$ and let $u_{\max} := \max_{1 \leq i \leq k-1}|u_i|.$
For fixed $k \geq 2$, we first claim that for $T := \frac{x}{(\log x)^{k+1}}$, the contribution of 
$k$-tuples that contain at least one element $a_n \leq T$ is negligible:

Indeed, we drop the roughness condition here completely, and find the upper bound
\[\frac{1}{\Phi(x,f(x))}\#\left\{(m, \bm{n}) \in [1,T]\times [1,x]^{k-1}: \lVert (m\mathds{1}_{k-1} - \bm{n})\alpha)\rVert_{\infty} \leq \frac{2u_{\max}}{\Phi(x,f(x))}\right\},\]
We can further upper-bound the above by
\[\frac{1}{\Phi(x,f(x))}\sum_{m \leq T} \left(\#\left\{n \leq x: \{n\alpha\} \in \left[m\alpha - \frac{2u_{\max}}{\Phi(x,f(x))}, m\alpha - \frac{2u_{\max}}{\Phi(x,f(x))}\right]\right\}\right)^{k-1}.\]
Since the minimal gap of $(\{n\alpha\})_{n \leq x}$ for badly approximable numbers is of size $\gg_{\alpha} 1/x$, we get

\[\#\left\{n \leq x: \{n\alpha\} \in \left[m\alpha - \frac{2u_{\max}}{\Phi(x,f(x))}, m\alpha - \frac{2u_{\max}}{\Phi(x,f(x))}\right]\right\}
\ll_{\alpha,\bm{R}} \frac{x}{\Phi(x,f(x))},\]
with the implied constant being uniform in $m$.
This proves that the contribution of tuples with at least one element $a_n \leq T$ is bounded from above by
\[\ll_{\alpha,\bm{R},k} \frac{T}{\Phi(x,f(x))}\left(\frac{x}{\Phi(x,f(x))}\right)^{k-1}.\]
Observe that $\Phi(x,f(x)) \sim x\prod_{p \leq f(x)}\left(1 - \frac{1}{p}\right) \gg \frac{x}{\log x}$, thus we have
with the choice of $T = \frac{x}{(\log x)^{k+1}}$ that the above is $o(1)$.
Thus $R_k$ is asymptotically equal to only counting the tuples
$(m_1,\ldots,m_k)$ with $\min_{1 \leq i \leq k}m_i > \frac{x}{(\log x)^{k+1}}.$

We now write $z^- = f(\frac{x}{(\log x)^{k+1}})$ and $z = f(x)$.
Since by the preceding discussion we can remove the contribution of elements $\leq \frac{x}{(\log x)^{k+1}}$, we can bound

\[\begin{split}&\frac{1}{\Phi(x,z^{-})}\#\left\{(m,\bm{n}) \in [1,x]^{k}_{\neq}: P^{-}((m,\bm{n})) > z, m\bm{1}_{k-1} - \bm{n} \in \frac{\bm{R}}{\Phi(x,z^-)}\cdot \frac{\Phi(x,z^-)}{\Phi(x,z)}\right\} + o(1)
\\&\leq
\frac{1}{N(x,f)}\#\left\{(m,\bm{n}) \in [1,x]^{k}_{\neq}: P^{-}(m) > f(m), P^{-}(n_i) > f(n_i), m\bm{1}_{k-1} - \bm{n} \in \frac{\bm{R}}{N(x,f)}\right\}
  = R_k(s,N)\\& \leq \frac{1}{\Phi(x,z)}\#\left\{(m,\bm{n}) \in [1,x]^{k}_{\neq}: P^{-}((m,\bm{n})) > z^-, m\bm{1}_{k-1} - \bm{n} \in \frac{\bm{R}}{\Phi(x,z^-)}\right\} + o(1),
 \end{split}\]
 where $N(x,f) = \#\{n \leq x: P^{-}(n) \geq f(n)\}$.
Observe that 
$\Phi(x,z) \leq N(x,f) \leq \Phi(x,z^{-})$ and 
$\Phi(x,z) \sim \Phi(x,z^{-})$ since
$\frac{\log z}{\log z^-} \to 1$.
 
 Applying Theorem \ref{main_thm} for triangular arrays (proven in Section \ref{sec_step_fct} above), once with $f(x)$ and once with $\tilde{f}(x) := f(x/ (\log x)^{k+1})$, we obtain

 \[\frac{1}{\Phi(x,z^{-})}\#\left\{(m,\bm{n}) \in [1,x]^{k}_{\neq}: P^{-}((m,\bm{n})) > z, m\bm{1}_{k-1} - \bm{n} \in \frac{\bm{R}}{\Phi(x,z^-)}\cdot \frac{\Phi(x,z^-)}{\Phi(x,z)}\right\} \sim \vol(\bm{R}),\]
 as well as
 \[\frac{1}{\Phi(x,z)}\#\left\{(m,\bm{n}) \in [1,x]^{k}_{\neq}: P^{-}((m,\bm{n})) > z, m\bm{1}_{k-1} - \bm{n} \in \frac{\bm{R}}{\Phi(x,z^-)}\right\} \sim \vol(\bm{R}).\]
This finally concludes the proof of Theorem \ref{main_thm}.

\section{Proof of Theorem \ref{dioph_ass_thm}}
For simplicity, we prove Theorem \ref{dioph_ass_thm} again for triangular arrays $(\{a_{n,x}\alpha\})_{n\leq x}$. The step to move to the sequence $(\{a_{n}\alpha\})_{n \in \mathbb{N}}$ instead can be established as in Section \ref{sec_arr_to_seq}. We adopt a simplified version of the strategy employed by Walker \cite{walker_primes}. By assumption \eqref{dioph_ass}, for any $\varepsilon > 0$, we have infinitely many integers $x$ such that
\begin{equation}\label{good_approx_N}\lVert x\alpha\rVert < \frac{\varepsilon}{x \log f(x))},\end{equation}
which trivially implies
\[\lVert h x\alpha\rVert < \frac{\varepsilon}{x}, \quad 1 \leq h \leq \log(f(x)).\]

Using $\Phi(x,f(x)) \asymp \frac{x}{\log f(x)}$, we obtain for sufficiently large $x$ (that satisfies \eqref{good_approx_N}) and suitably chosen constants $c_1,c_2,c_3,c_4 >0$ that

\[\begin{split}&\#\left\{1 \leq n\neq m \leq \Phi(x\log f(x),f(x\log f(x))): \lVert (a_n^{(f)} - a_{m}^{(f)})\alpha \rVert \leq \frac{\varepsilon}{\Phi(x\log f(x),f(x\log f(x)))}\right\} \\&\gg \#\left\{1 \leq n \neq m \leq c_1 x \log f(x): P^{-}(n) > f(x\log f(x)), P^{-}(m) > f(x\log f(x)), \lVert (n-m)\alpha \rVert \leq \frac{c_2\varepsilon}{x}\right\} 
\\&\gg \hspace{-1cm}\sum_{c_3\log(f(x)) \leq h \leq c_4\log(f(x))} \#\left\{1 \leq n,m \leq c_1x \log f(x): P^{-}(n) > f(x\log f(x)), P^{-}(m) > f(x\log f(x)), m-n = x h\right\}.
\end{split}
\]
Applying Lemma \ref{sieve_k_tuple} with $k = 2$ implies for even $h$
\[\#\left\{1 \leq n,m \leq c_1x \log f(x): P^{-}(n) > f(x\log f(x)), P^{-}(m) > f(x\log f(x)), m-n = x h\right\}
\gg \frac{x}{\log f(x)}.
\]
 Thus along the subsequence that satisfies \eqref{good_approx_N}, we get

\[\frac{\#\left\{1 \leq n\neq m \leq \Phi(x\log f(x),f(x\log f(x))): \lVert (a_n^{(f)} - a_{m}^{(f)})\alpha \rVert \leq \frac{\varepsilon}{\Phi(x\log f(x),f(x\log f(x)))}\right\}}{\Phi(x\log f(x),f(x\log f(x)))} \gg 1,
\]
with the implied constant independent of $\varepsilon$. Choosing $\varepsilon$ sufficiently small now contradicts the assumption that the above converges to $2 \varepsilon$,
and thus concludes the first part of Theorem \ref{dioph_ass_thm}. For the statement on the Lebesgue measure, we observe that \eqref{dioph_ass} holds for almost every $\alpha$: This fact can be recovered straightforwardly from Khintchine's famous theorem in metric Diophantine approximation \cite{Khi_thm}: It states that for $\psi: \N \to [0,\infty)$ a monotonically increasing function with $\sum_{n \in \mathbb{N}}\psi(n) = \infty,$ we have for almost every $\alpha$, $\lVert n \alpha\rVert < \psi(n)$ for infinitely many $n$. Since $\sum_{n \in \N}\frac{\varepsilon}{n \log f(n)} = \infty$ for all $\varepsilon > 0$, this concludes the proof.

% Requires \usepackage{float} for the [H] figure placement.
%\appendix
\section*{Appendix: Numerical observations on rough and prime rotations}
%\label{app:numerics}

\begin{figure}[H]
\centering
\includegraphics[width=.91\textwidth]{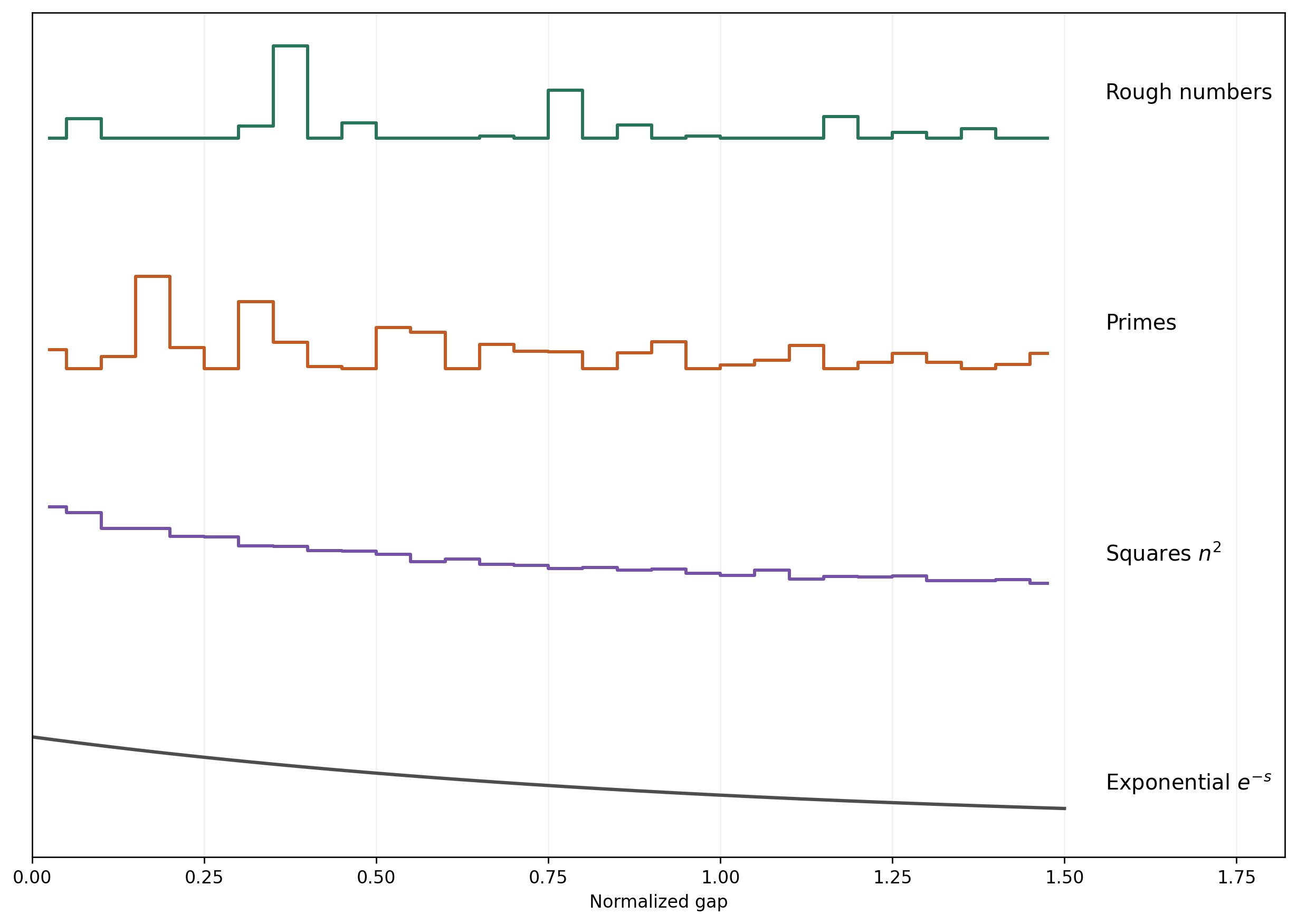}
\caption{Finite gap distributions at $M=10^8$ for the rough-number, prime, and square rotations, together with the exponential density $e^{-s}$. Each row is individually rescaled to unit height.}
\label{fig:unitheight}
\end{figure}

We complement Theorem~\ref{main_thm} with some numerical observations on the rough-number rotations considered in the theorem and on the prime rotation conjectured to have the same Poissonian behaviour. Throughout, we take
\[
    M=10^8,
    \qquad
    \alpha=\frac{1+\sqrt5}{2}=[1;1,1,\ldots],
\]
and for the rough numbers use the concrete choice
\[
    f(n)=\exp(\sqrt{\log n}),
\]
already mentioned after Theorem~\ref{main_thm}. Up to the cutoff $M$, this gives $13\,137\,988$ rough numbers, while there are $\pi(M)=5\,761\,455$ primes. We also include the square rotation $\{n^2\alpha\}$, whose conjectured Poissonian fine-scale behaviour was discussed in the Introduction in connection with the Conjectures made by Rudnick--Sarnak--Zaharescu \cite{rud_sar,rhz} and Heath--Brown \cite{hb}. Figure~\ref{fig:unitheight} displays the normalized nearest-neighbour gap distributions of these three rotations. Each row is rescaled to unit height and vertically shifted, so that the figure emphasizes the shape of the distributions at the finite cutoff.

At the cutoff $M=10^8$, the fine-bin gap histograms in Figure~\ref{fig:unitheight} still show pronounced discrete structure, and in particular are not yet visually close to the limiting exponential law. Before explaining this phenomenon more closely

We next turn to a more arithmetic finite-cutoff comparison for the prime rotation, and later return to explain the origin of the structure visible in the gap histograms: We compare the prime correlations with their finite Hardy--Littlewood predictions. This comparison is natural from the structure of the proof: Lemma~\ref{sieve_k_tuple} supplies the uniform $k$-tuple asymptotic for the rough numbers, while the corresponding prime input would require a sufficiently uniform form of the Hardy--Littlewood $k$-tuple conjecture. This is the same type of uniform Hardy--Littlewood input appearing in Gallagher's conditional derivation of Poissonian statistics for the primes on the real line \cite{Gall_primes}. Thus the numerical comparison below may be viewed as a finite-cutoff test of the prime analogue of the sieve input in Lemma~\ref{sieve_k_tuple}. Once the sieve input in circumvented, the proof of Poissonian behaviour for the primes would follow along the lines of Theorem \ref{main_thm} since the respective singular series coincide,
and the averaging technique of the Diophantine Bohr set becomes only easier.

For the asymptotic discussion in the main part of the paper, the prime density $1/\log M$ is sufficient. For the finite numerical comparison, we retain the variation of the logarithmic density over $[2,M]$, which gives a more accurate prediction at the cutoff considered here. Using the notation $g_{\bm h}(p)$, $h^+_{>0}$, and $h^-_{<0}$ introduced above, for $(0,\bm h)\in\Z^k_{\neq}$ we set
\[
 \operatorname{HL}_M(\bm h)
 :=
 \left(\prod_p\frac{1-g_{\bm h}(p)/p}{(1-1/p)^k}\right)
 \int_{2-h^-_{<0}}^{M-h^+_{>0}}
 \frac{\mathrm{d}t}{\log t\prod_{j=1}^{k-1}\log(t+h_j)},
\]
with value $0$ if the integration interval is empty. Thus $\operatorname{HL}_M(\bm h)$ is the usual finite Hardy--Littlewood main term for the number of prime tuples $p,p+h_1,\ldots,p+h_{k-1}$ contained in $[2,M]$.
For $N=\pi(M)$, define
\[
 R_k^{\rm HL}([-s,s]^{k-1},N)
 :=\frac1N
 \sum_{\substack{
      \bm h\in[-M,M]^{k-1}\\
      (0,\bm h)\in\Z^k_{\neq}\\
      \lVert h_j\alpha\rVert\leq s/N\; (1\leq j\leq k-1)}}
 \operatorname{HL}_M(\bm h).
\]
In other words, $R_k^{\rm HL}$ is obtained from the difference representation of $R_k([-s,s]^{k-1},N)$ by replacing the exact prime $k$-tuple count for each difference vector $\bm h$ by its finite Hardy--Littlewood prediction. The singular-series averaging established in the paper implies, for fixed $k$ and $s$,
\[
    R_k^{\rm HL}([-s,s]^{k-1},N)\longrightarrow(2s)^{k-1}.
\]

Thus it remains to show that $R_k([-s,s]^{k-1},N) \sim R_k^{\rm HL}([-s,s]^{k-1},N)$ which would be a consequence of the uniform Hardy--Littlewood Conjecture described above, but we also provide numerical support: Indeed, with $M = 10^8$ and $\alpha = [1;1,1,\ldots]$, the respective values are as follows:

\begin{center}
\setlength{\tabcolsep}{8pt}
\begin{tabular}{c|c|c|c|c}
$k$ & $s$ & $R_k([-s,s]^{k-1},N)$ & $R_k^{\rm HL}([-s,s]^{k-1},N)$ & $(2s)^{k-1}$\\
\hline
$2$ & $0.5$ & $0.854003$ & $0.853988$ & $1$\\
    & $1$   & $1.937222$ & $1.937291$ & $2$\\
    & $2$   & $3.701886$ & $3.701936$ & $4$\\
    & $3$   & $5.935788$ & $5.936450$ & $6$\\
    & $5$   & $9.873193$ & $9.874589$ & $10$\\
\hline
$3$ & $0.5$ & $0.637084$ & $0.636127$ & $1$\\
    & $1$   & $3.372202$ & $3.369405$ & $4$\\
    & $2$   & $13.066014$ & $13.060035$ & $16$\\
    & $3$   & $34.147779$ & $34.144958$ & $36$\\
    & $5$   & $95.420110$ & $95.425003$ & $100$\\
\hline
$4$ & $0.5$ & $0.400420$ & $0.399242$ & $1$\\
    & $1$   & $5.176273$ & $5.163474$ & $8$\\
    & $2$   & $43.703599$ & $43.655932$ & $64$\\
    & $3$   & $189.972642$ & $189.881178$ & $216$\\
    & $5$   & $902.295405$ & $902.142111$ & $1000$.
\end{tabular}
\end{center}
At all displayed points the relative difference between $R_k$ and $R_k^{\rm HL}$ is below $3 \cdot 10^{-4}$, whereas their common finite-cutoff value can still differ substantially from the Poissonian limit $(2s)^{k-1}$. This allows to isolate the unproven part of the conjecture from the dominating finite cutoff effect, providing numerical support for the Poissonian prime rotation conjecture across several correlation orders.

% \begin{figure}[H]
% \centering
% \includegraphics[width=.99\textwidth]{prime_HL_correlations_k234_v5.png}
% \caption{Prime correlations and their finite Hardy--Littlewood predictions at $M=10^8$ for $k=2,3,4$, normalized by the Poissonian value $(2s)^{k-1}$. The solid curves show $R_k^{\rm HL}$ evaluated on a fine mesh in $s$, the circles show the directly computed prime statistics at $s=0.5,1,2,3,5$, and the dashed line is the Poissonian value.}
% \label{fig:hlcomparison}
% \end{figure}
% Figure~\ref{fig:hlcomparison} shows close agreement between the prime statistics and their finite Hardy--Littlewood predictions for all three correlation orders.  The numerical values at the marked points, together with the corresponding Poissonian values, are

\begin{figure}[H]
\centering
\includegraphics[width=.93\textwidth]{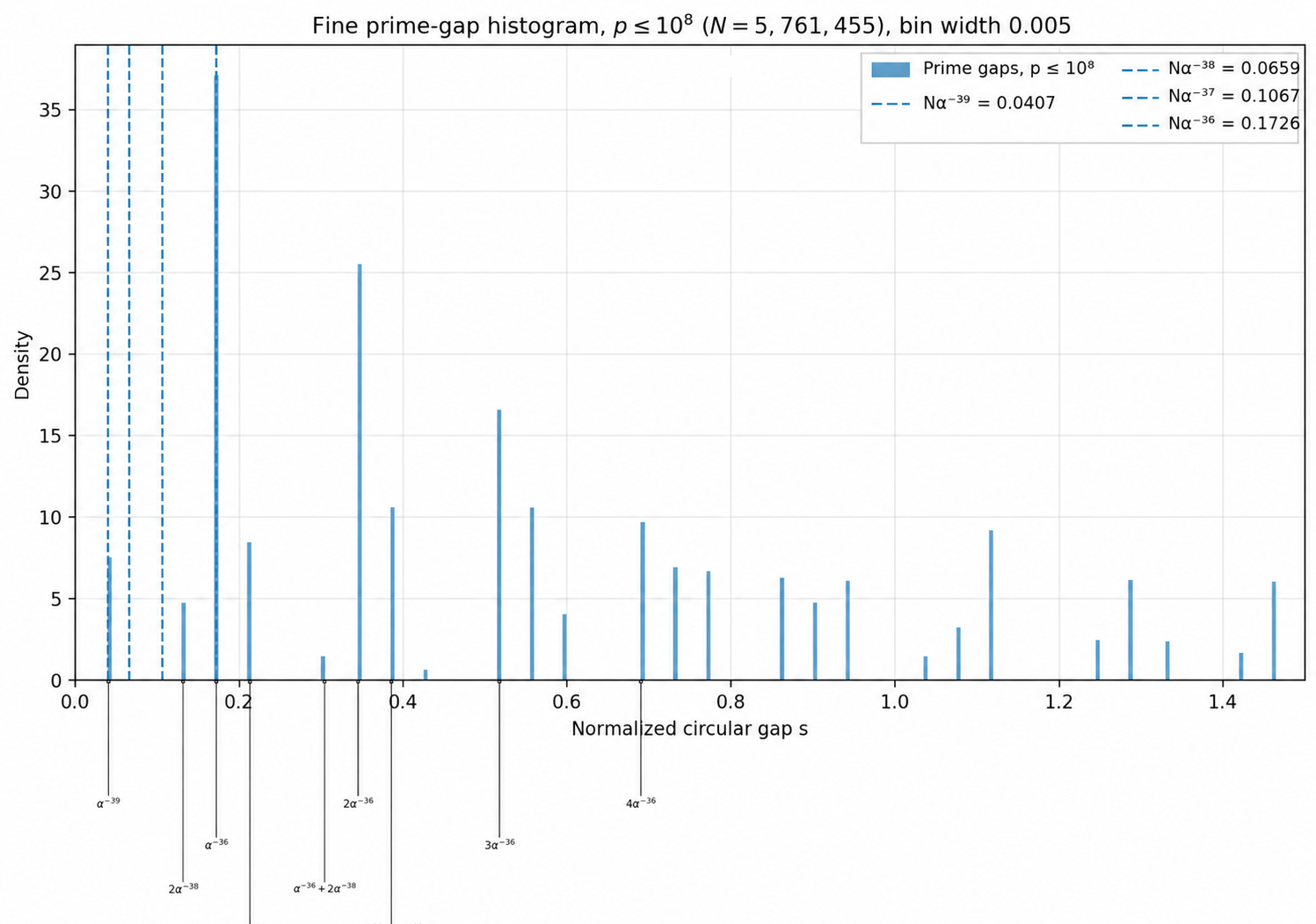}
\vspace{1mm}

\includegraphics[width=.93\textwidth]{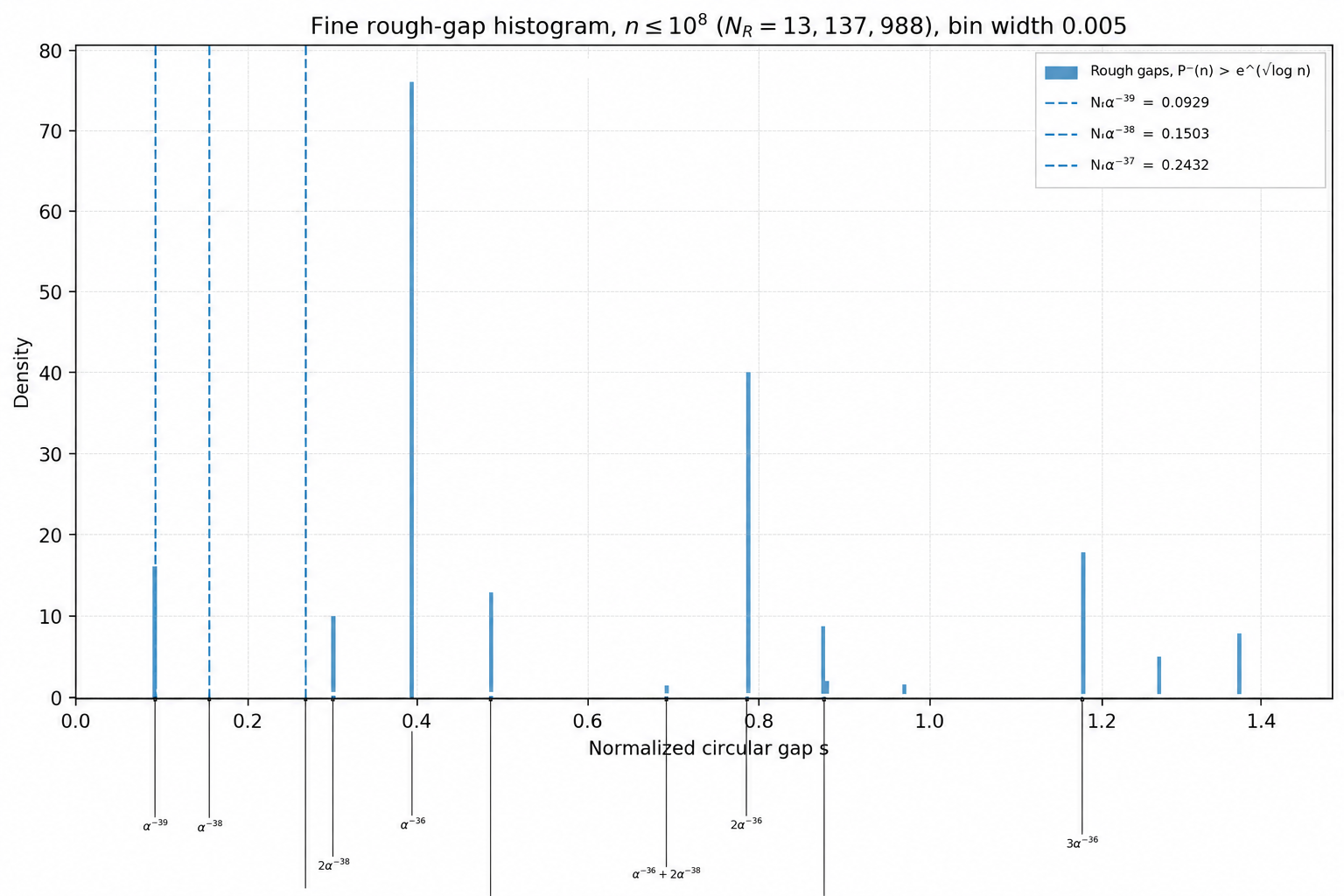}
\caption{Very fine gap histograms at $M=10^8$ Top: primes. Bottom: rough numbers. The dominant arithmetic family generated by $\alpha^{-36}$ and several secondary combinations of the elementary Three Gap lengths are marked.}
\label{fig:finegaps}
\end{figure}

We finally explain the visible structure in Figure~\ref{fig:unitheight}. For the Golden Ratio, the convergent denominators are the Fibonacci numbers and
\[
    F_{39}=63\,245\,986<M<F_{40}=102\,334\,155,
\]
so the full orbit $\{n\alpha:1\leq n\leq M\}$ has the three gap lengths
\[
    \alpha^{-39},\quad \alpha^{-38},\quad \alpha^{-37}.
\]
After restriction to a subsequence with $N$ retained points and normalization of the gaps by $N$, these elementary lengths are all of order $N/M$ --- which holds in general for $\alpha \in BAD$. 
Local arithmetic restrictions amplify the effect of these elementary lengths causing spikes: they weight particular combinations of the three lengths. In particular $\alpha^{-38}+\alpha^{-37}=\alpha^{-36}$
and the corresponding Fibonacci displacement $F_{36}$ is even, a requirement for prime differences (except the negligible $p = 2)$. This produces the especially pronounced family of spikes at multiples of $\alpha^{-36}$ visible in Figure~\ref{fig:finegaps}, together with secondary families arising from other admissible combinations.

For the three index sets used in Figure~\ref{fig:unitheight}, the normalized scales are
\[
\begin{array}{c|c|c}
\text{index set} & \text{asymptotic normalized mesh scale} & M=10^8\\ \hline
\{n\leq M:P^-(n)>f(n)\}
& \displaystyle \frac{N(M,f)}{M}\sim\frac{e^{-\gamma}}{\sqrt{\log M}}
&0.1314\\[2mm]
\{p\leq M:p\text{ prime}\}
& \displaystyle \frac{\pi(M)}{M}\sim\frac1{\log M}
&0.0576\\[2mm]
\{n^2\leq M\}
& \displaystyle \frac{\lfloor\sqrt M\rfloor}{M}\sim\frac1{\sqrt M}
&0.0001.
\end{array}
\]

\begin{figure}[H]
\centering
\includegraphics[width=.94\textwidth]{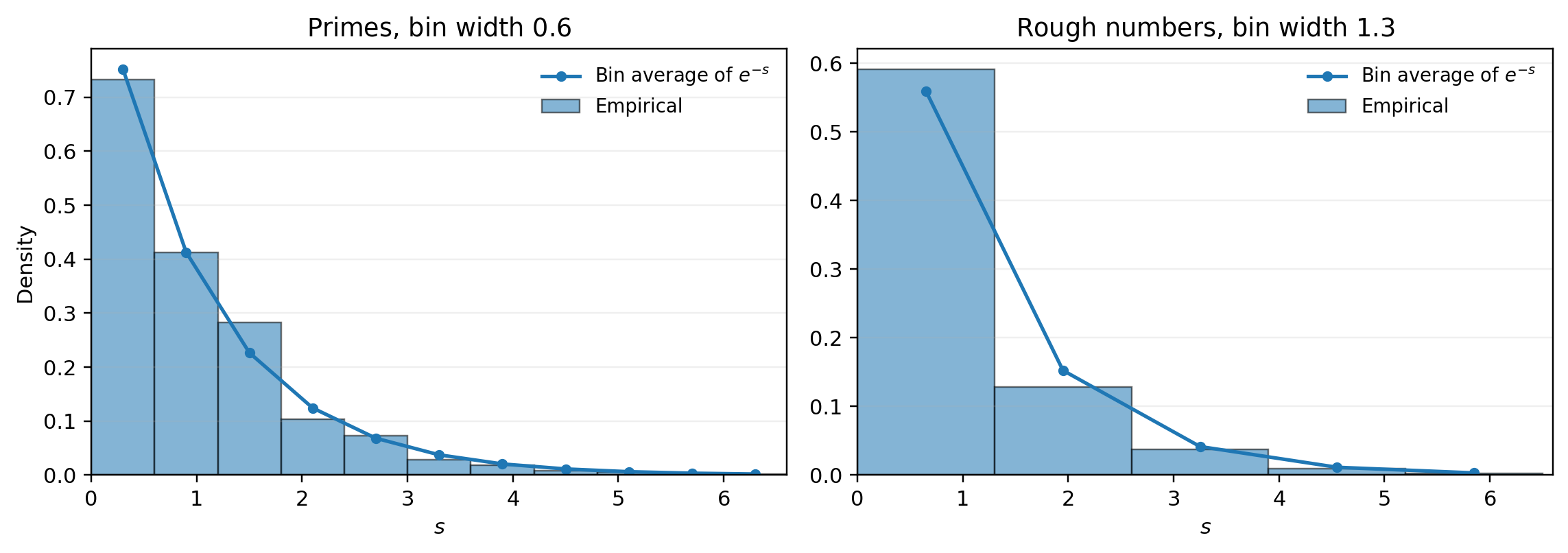}
\caption{The prime and rough-number gap distributions at $M=10^8$ after coarse-graining. The bin widths $0.6$ for the primes and $1.3$ for the rough numbers are approximately ten times the corresponding normalized mesh scales $N/M$. The line gives the average of $e^{-s}$ over each bin.}
\label{fig:coarsegaps}
\end{figure}

Thus the normalized mesh is already very fine for the square rotation at $M=10^8$, whereas for the rough-number and prime rotations it remains visible at the same cutoff. Figure~\ref{fig:coarsegaps} shows the same prime and rough-number gap data after coarse-graining at bin widths chosen to be approximately ten times their respective normalized mesh scales. At this coarser resolution, the exponential envelope is much more clearly visible.

While these bins look unusually big here, note that for the prime and rough-number sequences, any fixed multiple of the mesh scale still tends to $0$ as $M\to\infty$. This contrasts with integer sequences of positive density, for which $N/M$ stays away from zero, and thus the normalized Three Gap structure remains even asymptotically visible. The convergence of $N/M$ to $0$ for primes and rough numbers is, however, very slow. 

% For a further numerical summary of the finite gap distributions, let $F_M$ denote the empirical distribution function of the normalized nearest-neighbour gaps and define the Kolmogorov--Smirnov distance from the exponential distribution by
% \[
%     D_{\rm KS}:=\sup_{s\geq0}\left|F_M(s)-(1-e^{-s})\right|.
% \]
% At $M=10^8$ we obtain
% \[
% \begin{array}{c|c}
% \text{index set} & D_{\rm KS}\\ \hline
% \{n\leq M:P^-(n)>f(n)\} & 0.22133\\
% \{p\leq M:p\text{ prime}\} & 0.11208\\
% \{n^2\leq M\} & 0.03141.
% \end{array}
% \]
% These values are consistent with the markedly different mesh scales above, and we interpret the sizeable finite-$M$ deviations for the prime and rough-number rotations in this context. 
%The slow decay of the mesh is particularly relevant for fine-bin histograms: 
At a bin size of order $10^{-2}$, pronounced finite-size effects should remain visible as long as the normalized mesh is comparable with this scale. Already the relation $1/\log M=10^{-2}$ for the primes corresponds to $M=e^{100}\approx 3\cdot10^{43}$ while the analogous cutoff for the rough-number rotation is of order $10^{1369}$. Thus reaching a regime in which the mesh is substantially finer than the $10^{-2}$ plotting scale would require cutoffs far beyond reasonable direct computation. The finite Hardy--Littlewood comparison of $R_k$ and $R_k^{HL}$ therefore provides the more meaningful numerical stress-test of the Poissonian Prime rotation Conjecture.

\bibliographystyle{abbrv}
\bibliography{bibliography.bib}
\Addresses

\end{document}